\documentclass[twoside,12pt]{article}
\usepackage[]{amsmath,amsthm,amssymb}
\usepackage[latin1]{inputenc}

\begin{document}

\title{{\Large\bf On the Macdonald-type function and its relation with index transforms and orthogonal polynomials}}

\author{Semyon  YAKUBOVICH}
\maketitle

\markboth{\rm \centerline{ Semyon   YAKUBOVICH}}{}
\markright{\rm \centerline{Macdonald-Type Function}}

\begin{abstract} {We continue to investigate properties of the function $M_\nu(z)$ which is associated  with the Macdonald function $K_\nu(z)$ in terms of the corresponding Fourier integral.  In particular, recurrence relations for this function and its derivatives are obtained, involving properties of the associated Laguerre polynomials. Multiple orthogonal polynomials related to the scaled Macdonald-type  weights $ \hat{\rho}_{\nu}(x)= 2 x^{\nu/2} M_\nu\left(2\sqrt x\right),  x >0$  are investigated. }

\end{abstract}
\vspace{4mm}

{\bf Keywords}: {\it Modified Bessel functions, Macdonald function, Exponential integral,
Kontorovich-Lebedev transform,  Fourier transform, Orthogonal polynomials }

{\bf AMS subject classification}:    44A15,  33B20, 33C10, 33C52

\vspace{4mm}

\section {Introduction and preliminary results}

As it is known [4],  the modified Bessel function or Macdonald function $K_\nu(x)$ can be defined by the integral representation

$$K_{\nu}(x)= \int_{0}^\infty e^{-x\cosh u} \cosh(\nu u) du,\quad x >0, \nu \in \mathbb{C}.\eqno(1.1)$$
The pure imaginary $\nu=i\tau,\ \tau \in \mathbb{R}$ corresponds the kernel $K_{i\tau}(x)$ of the reciprocal Kontorovich-Lebedev transforms 

$$F(\tau)=\int_0^\infty K_{i\tau}(x) f(x) {dx\over x},\eqno(1.2)$$

$$f(x)= {2\over \pi^2} \int_0^\infty \tau\sinh(\pi\tau) K_{i\tau}(x) F(\tau) d\tau\eqno(1.3)$$
which generate the class of the index transforms [8], where the integration is realized with respect to the argument and the index (a parameter) of the kernel function. The function $K_\nu(x)$ is a fundamental solution of the Bessel homogeneous second order differential equation

$$x^2 y^{\prime\prime} (x) + x y^{\prime} (x) - (x^2+\nu^2) y(x)= 0.\eqno(1.4)$$
It means that $K_\nu(x)$ is the eigenfunction of the differential operator

 $${\cal A} \equiv x^2- x{d\over dx} x {d\over dx},\eqno(1.5)$$
i.e. this yields
$${\cal A} \ K_{\nu}(x)= - \nu^2 K_{\nu}(x).\eqno(1.6)$$
Moreover, $K_{\nu}(x)$  remains bounded as $x$ tends to
infinity on the real positive line. It has the asymptotic behaviour [4]
$$ K_\nu(x) = \left( \frac{\pi}{2x} \right)^{1/2} e^{-x} [1+
O(1/x)], \qquad x \to +\infty,\eqno(1.7)$$
and near the origin
$$ K_\nu(x) = O\left ( x^{-|{\rm Re}\nu|}\right), \ x \to 0,\eqno(1.8)$$
$$K_0(x) = -\log x + O(1), \ x \to 0. \eqno(1.9)$$
Our main goal in the sequel is to focus on the properties of the companion function $M_\nu(x)$ which was introduced and studied in [3] via the integral representation

$$M_{\nu}(x)= \int_{0}^\infty e^{-x\cosh u} \sinh(\nu u) du,\quad x >0, \nu \in \mathbb{C}.\eqno(1.10)$$
This function is a solution of the following non-homogeneous Bessel differential equation
$$x^2 y^{\prime\prime} (x) + x y^{\prime} (x) - (x^2+\nu^2) y(x)= \nu e^{-x}.\eqno(1.11)$$
Furthermore,  it obeys the following recurrence and differential recurrence relations (see [3])

$$x \left[ M_{\nu+1}(x) - M_{\nu-1}(x) \right] = 2 \left[ e^{-x} + \nu M_\nu(x)\right],\eqno(1.12)$$

$$ M_{\nu+1}(x) + M_{\nu-1}(x) = - 2 {dM_\nu\over dx},\eqno(1.13)$$

$$x {dM_\nu\over dx} +  \nu M_\nu(x) = - e^{-x} - x M_{\nu-1}(x),\eqno(1.14)$$

$$x {dM_\nu\over dx} -  \nu M_\nu(x) = e^{-x} - x M_{\nu+1}(x),\eqno(1.15)$$

$${d\over dx} \left[ x^\nu M_\nu (x)\right] = - x^\nu M_{\nu-1} (x)  - x^{\nu-1} e^{-x},\eqno(1.16)$$

$${d\over dx} \left[ x^{-\nu} M_\nu (x)\right] = - x^{-\nu} M_{\nu+1} (x)  + x^{-\nu-1} e^{-x}.\eqno(1.17)$$
The pure imaginary case of the index $\nu=i\tau$ produces the function $M_{i\tau}$ which is the kernel of following reciprocal pairs of the index transforms, involving the modified Bessel function of the first kind $I_\nu(x)$ [4] 

$$g(x)= \int_0^\infty \tau M_{i\tau}(x) f(\tau)d\tau,\eqno(1.18)$$

$$f(\tau)= - {1\over \pi}  \int_0^\infty \bigg[ I_{i\tau}(x) +  I_{-i\tau}(x)\bigg] g(x) {dx\over x},\eqno(1.19)$$

$$\Phi(\tau) = \int_0^\infty M_{i\tau} (x) \varphi(x) {dx\over x},\eqno(1.20)$$

$$\varphi(x)= - {1\over \pi}  \int_0^\infty \tau \bigg[ I_{i\tau}(x) +  I_{-i\tau}(x)\bigg] \Phi(\tau) d\tau,\eqno(1.21)$$
which were introduced and inverted in [6] (slightly corrected) and [9], respectively, in suitable spaces of functions.  To end this section  it is worth to mention the sequence of polynomials $\{p_n(x)\}_{n\ge 0}$ related to the Kontorovich-Lebedev transform (1.2)  (cf. [10]) which is defined in terms of the $n$th iteration of the differential operator (1.5)
$$p_n(x)= (-1)^n e^{x}{\mathcal A}^n\  e^{-x}, \ n \in \mathbb{N}_0.\eqno(1.22)$$
Evidently, $p_0(x)=1, p_1(x)= -x$. It has the integral representation
$$p_n(x)= {2(-1)^n\over \pi} e^{x}\int_{0}^\infty \tau^{2n} K_{i\tau}(x)\ d\tau,\ x >0,\eqno(1.23) $$
 and satisfies  the differential recurrence relation of the form
$$p_{n+1}(x)= x^2p_n^{\prime\prime}(x) + x(1-2x)p_n^\prime(x)- xp_n(x), \ n= 0,1,2,\dots \ .\eqno(1.24)$$
The generating function for this system of polynomials is given by the series
$$e^{- 2x\sinh^2( y/2)} = \sum_{n=0}^\infty {p_n(x)\over (2n)!}\ y^{2n}\eqno(1.25)$$
with a positive convergence radius.  Hence

$$e^{- x\cosh y} = e^{-x} \sum_{k=0}^\infty {p_k(x)\over (2k)!}\ y^{2k},$$
and via termwise differentiation by $y$ within the convergence radius, we derive

$$D_y^{2n} \bigg[ e^{-x\cosh y}\bigg] = e^{-x} \sum_{k=0}^\infty {p_{k+n}(x)\over (2k)!}\ y^{2k}.\eqno(1.26)$$
Polynomials  $p_n$ obey the  recurrence relation 

$$p_{n+1}(x)= -x \sum_{k=0}^n \binom{2n+1}{2k} p_k(x),\quad n \in \mathbb{N}_0\eqno(1.27)$$
and for the derivatives  $p^\prime_n$  we have 

$$p^\prime_{n}(x)= - \sum_{k=0}^{n-1} \binom{2n}{2k} p_k(x),\quad n \in \mathbb{N}.$$

\section{The use of the associated Laguerre polynomials}

Returning to (1.11), we write the differential equation for the function $M_\nu(x)$ in the operator form

$${\cal A}  M_{\nu}(x)= - \nu^2 M_{\nu}(x)- \nu e^{-x}.\eqno(2.1)$$
Hence, employing (1.22), we get immediately from (2.1)

$$ {\cal A}^2  M_{\nu}(x)=  \nu^4  M_{\nu}(x) + \nu^3 e^{-x} + \nu e^{-x} p_1(x).\eqno(2.2)$$
Consequently, for the $n$th iteration of the operator ${\cal A}$ (${\cal A}^0 \equiv I$)  we derive the formula

$$ (-1)^n  {\cal A}^n  M_{\nu}(x)=   \nu^{2n}   M_{\nu}(x)  +  e^{-x}  \sum_{k=1}^{n}  \nu^{2k-1} p_{n-k}(x),\quad n \in \mathbb{N}_0.\eqno(2.3)$$
Meanwhile, identity (1.16) permits an involvement of the associated Laguerre polynomials [4] $L_n^\nu(x),\ n \in \mathbb{N}_0$. Indeed, using the Rodrigues formula 

$$L_n^\nu(x) = {e^{x} x^{-\nu}\over n!} {d^n\over dx^n} \left[ e^{-x} x^{n+\nu}\right],\eqno(2.4)$$
identity (1.16) yields

$$ {d^n\over dx^n} \left[ x^{n+1} {d\over dx} \left[ x^\nu M_\nu (x)\right] \right] + {d^n\over dx^n} \left[ x^{\nu+n+1} M_{\nu-1} (x) \right] =  - n! e^{-x} x^{\nu} L_n^\nu(x).\eqno(2.5)$$
Denoting by 

$$\omega^\alpha_{n,\nu} (x) = {d^n\over dx^n} \left[ x^{n+\alpha} M_\nu (x)\right],\eqno(2.6)$$
and appealing to (1.13), we write the left-hand side  of (2.5) in the form

$$  {d^n\over dx^n} \left[ x^{n+1} {d\over dx} \left[ x^\nu M_\nu (x)\right] \right] + {d^n\over dx^n} \left[ x^{\nu+n+1} M_{\nu-1} (x) \right] $$

$$=  {d^{n+1}\over dx^{n+1}} \left[ x^{\nu+n+1} M_\nu (x)\right] - (n+1)  {d^{n}\over dx^{n}} \left[ x^{\nu+n}  M_\nu(x) \right]  - 2 {d^{n+1}\over dx^{n+1}} \left[ x^{\nu+n+1}  M_\nu(x) \right]$$

$$+2 (n+1+\nu) {d^{n}\over dx^{n}} \left[ x^{\nu+n}  M_\nu(x) \right]   - {d^n\over dx^n} \bigg[ x^{\nu+n+1}  M_{\nu+1}(x)\bigg] $$

$$= - \omega^\nu_{n+1,\nu} (x) - \omega^\nu_{n,\nu+1} (x) +  (n+1+ 2\nu) \omega^\nu_{n,\nu} (x).$$
Hence identity (2.5) can be written in the final form

$$ \omega^\nu_{n+1,\nu} (x) + \omega^\nu_{n,\nu+1} (x) -  (n+1+ 2\nu) \omega^\nu_{n,\nu} (x) =  n! e^{-x} x^{\nu} L_n^\nu(x).\eqno(2.7)$$
However, the three-term recurrence relation for the associated Laguerre polynomials [4]

$$x L_n^\nu(x) = (2n+\nu+1) L_n^\nu(x) - (n+1) L_{n+1}^\nu(x) - (n+\nu) L_{n-1}^\nu(x),\eqno(2.8)$$
permits to rewrite (2.7) in the form

$$ x \bigg[ \omega^\nu_{n+1,\nu} (x) + \omega^\nu_{n,\nu+1} (x) -  (n+1+ 2\nu) \omega^\nu_{n,\nu} (x) \bigg] $$

$$= (2n+\nu+1) \bigg[ \omega^\nu_{n+1,\nu} (x) + \omega^\nu_{n,\nu+1} (x) -  (n+1+ 2\nu) \omega^\nu_{n,\nu} (x) \bigg]$$

$$-  \omega^\nu_{n+2,\nu} (x) - \omega^\nu_{n+1,\nu+1} (x) +  (n+2+ 2\nu) \omega^\nu_{n+1,\nu} (x)$$

$$- n(n+\nu) \bigg[ \omega^\nu_{n,\nu} (x) + \omega^\nu_{n-1,\nu+1} (x) -  (n+ 2\nu) \omega^\nu_{n-1,\nu} (x) \bigg],$$
i.e. after simplification we deduce the identity

$$  \omega^\nu_{n+2,\nu} (x) + \omega^\nu_{n+1,\nu+1} (x) $$

$$ =  \bigg[  (n+1+ 2\nu) ( x  - 2n-\nu-1)  - n(n+\nu) \bigg]  \omega^\nu_{n,\nu} (x) $$

$$ + \left[ 3 (n+1+ \nu) -x\right] \omega^\nu_{n+1,\nu} (x)+ (2n+\nu+1-x) \omega^\nu_{n,\nu+1} (x) $$

$$- n(n+\nu) \bigg[  \omega^\nu_{n-1,\nu+1} (x) -  (n+ 2\nu) \omega^\nu_{n-1,\nu} (x) \bigg].$$
Further,  employing the Mellin-Barnes integral representation for the function  $ e^x M_\nu (x)$ [6]

$$M_\nu(x) =  {e^{-x} \sin(\pi\nu)\over 2\pi^{3/2} i} \int_{\gamma-i\infty}^{\gamma+i\infty} {\Gamma (s+\nu)\Gamma(s-\nu) \Gamma(s)\Gamma(1-s)\over \Gamma(s+1/2)} (2x)^{-s} ds,\eqno(2.9)$$
where $\ x > 0,\  |{\rm Re} \nu | < \gamma < 1$, $\Gamma(z)$ is the Euler gamma function [4], we appeal to the generalized Mellin-Parseval equality [6] and Entry 8.4.23.27 in [5, Vol. III] to establish the integral representation for $M_\nu(x)$ in terms of the square of Macdonald's  function

$$M_\nu(x) =  {2\over \pi} \  e^{-x} \sin(\pi\nu) \int_0^\infty  e ^{-t} K^2_{\nu}\left(\sqrt{2x t}\right) dt,\ x >0,\  |{\rm Re} \nu |  < 1. \eqno(2.10)$$
It can be written in terms of the square of the scaled Macdonald function [11] $\rho_\nu(x)= 2 x^{\nu/2} K_\nu(2\sqrt x)$. Precisely, we have

$$ x^\nu M_\nu(2x) =  { e^{-2x} \over 2\pi} \  \sin(\pi\nu) \int_0^\infty  e ^{-t} \rho^2_\nu(xt)\ t^{-\nu} dt. \eqno(2.11)$$
Integrating by parts, eliminating the integrated terms via the asymptotic behavior of the Macdonald function at infinity and near the origin [3] and since the derivative of the scaled Macdonald function has the form $\rho^\prime_\nu(x) = - \rho_{\nu-1} (x)$ [11], we deduce from (2.11)

$$ x^\nu M_\nu(2x) =  { e^{-2x} \over 2\pi} \  \sin(\pi\nu) \bigg[ {1\over 1-\nu} \int_0^\infty  e ^{-t} \rho^2_\nu(xt)\ t^{1-\nu} dt\bigg.$$

$$+ {2x \over 1-\nu} \int_0^\infty  e ^{-t} \rho_\nu(xt) \rho_{\nu-1} (xt) \ t^{1-\nu} \bigg] dt$$

$$=  { e^{-2x} \over 2\pi} \  \sin(\pi\nu) \bigg[ {1\over 1-\nu} \int_0^\infty  e ^{-t} \rho^2_\nu(xt)\ t^{1-\nu} dt- {x \over 1-\nu}  {d\over dx} \int_0^\infty  e ^{-t} \rho^2_\nu(xt)  \ t^{-\nu} \bigg] dt,$$
where the differentiation with respect to $x$ under the integral sign is possible for $x \ge x_0 > 0$ due to the absolute and uniform convergence.  Therefore, recalling (2.11), it gives the equality

$$  x^\nu \bigg[ (2x+1) M_\nu(2x) + x {d\over dx}  \left[M_\nu(2x) \right] \bigg]= { e^{-2x} \over 2\pi} \  \sin(\pi\nu) \int_0^\infty  e ^{-t} \rho^2_\nu(xt)\ t^{1-\nu} dt.\eqno(2.12)$$
Now, taking in mind the recurrence relation for the scaled Macdonald functions [11]

$$\rho_{\nu+1}(x)= \nu\rho_\nu(x) + x \rho_{\nu-1}(x),\eqno(2.13)$$
we write from (2.11)  

$$ \nu^2 x^{\nu} M_{\nu} (2x) =   { e^{-2x}   \sin(\pi\nu) \over 2\pi}  \int_0^\infty  e ^{-t} \bigg[ \rho^2_{\nu+1}(xt)+ 2x\rho_{\nu+1} (xt)  {d\over dx} \left[ \rho_\nu(xt)\right] + (xt)^2  \rho^2_{\nu-1}(xt) \bigg] \ t^{-\nu} dt .$$
Hence, employing (2.12), the latter equality will take the form

$$ \nu^2 x^{\nu} M_{\nu} (2x) =  - x^{\nu+2}  {d\over dx}  \bigg[ M_{\nu+1} (2x) + x M_{\nu-1} (2x)  \bigg] $$

$$ -   x^{\nu+1} \bigg[ (2x+1) M_{\nu+1} (2x) + 2  x^2 M_{\nu-1} (2x)\bigg] $$

$$+ { e^{-2x}  x  \sin(\pi\nu) \over \pi}  \int_0^\infty  e ^{-t}  \rho_{\nu+1} (xt)  {d\over dx} \left[ \rho_\nu(xt)\right]  \ t^{-\nu} dt .\eqno(2.14)$$
Meanwhile,

$$ { e^{-2x}  x  \sin(\pi\nu) \over \pi}  \int_0^\infty  e ^{-t}  \rho_{\nu+1} (xt)  {d\over dx} \left[ \rho_\nu(xt)\right]  \ t^{-\nu} dt $$

$$=    { e^{-2x}  x\ \sin(\pi\nu) \over 2\pi} \bigg[ -  {d^2\over dx^2}   \int_0^\infty  e ^{-t}  \rho^2_{\nu+1} (xt)   \ t^{-1-\nu} dt + 2  \int_0^\infty  e ^{-t}  \rho^2_\nu(xt)  \ t^{1-\nu} dt \bigg].$$
Consequently,  with the use of (2.11), (2.12), equality (2.14) becomes

$$ \nu^2  M_{\nu} (2x) =  - x^{2}  {d\over dx}  \bigg[ M_{\nu+1} (2x) + x M_{\nu-1} (2x)  \bigg] $$

$$ -   x \bigg[ (2x+1) M_{\nu+1} (2x) + 2  x^2 M_{\nu-1} (2x)\bigg] $$

$$+  2 x \bigg[ (2x+1) M_\nu(2x) + x {d\over dx}  \left[M_\nu(2x) \right] \bigg] +  e^{-2x}  x^{1-\nu}  {d^2\over dx^2} \bigg[  e^{2x}  x^{\nu+1}  M_{\nu+1} (2x)  \bigg].\eqno(2.15)$$
But, in turn,

$$e^{-2x}  x^{1-\nu}  {d^2\over dx^2} \bigg[  e^{2x}  x^{\nu+1}  M_{\nu+1} (2x)  \bigg] = e^{-2x}  x^{1-\nu}  {d\over dx} \bigg[  2 e^{2x}  x^{\nu+1}  M_{\nu+1} (2x) $$

$$+ (\nu+1) e^{2x}  x^{\nu}  M_{\nu+1} (2x) +e^{2x}  x^{\nu+1} {d\over dx} [ M_{\nu+1} (2x) ] \bigg] $$

$$=  \left(\nu(\nu+1)+ 4   x (\nu+x+1) \right) M_{\nu+1} (2x) + 2x( 2  x+\nu+1)  {d\over dx} [ M_{\nu+1} (2x) ]$$

$$  +  x^{2} {d^2\over dx^2} [ M_{\nu+1} (2x) ].$$
Therefore,  invoking (1.14),  we establish from (2.15) after simplification the following differential recurrence relation
$$ x^{2} {d^2\over dx^2} [ M_{\nu+1} (2x) ]+   \bigg(3  x^2 + 2x(\nu+1) \bigg) {d\over dx} [ M_{\nu+1} (2x) ]+  \bigg(2x^2+ (4\nu+3) x+ \nu(\nu+1) \bigg) M_{\nu+1} (2x) $$

$$= - {x^{3}\over 2} {d^2\over dx^2}  [M_\nu(2x) ]  - \bigg( x+ {5+\nu\over 2} \bigg) x^2 {d\over dx}  \left[M_\nu(2x) \right] - \bigg(  (\nu+4) x^2+ 2x - \nu^2 \bigg)  M_\nu(2x).\eqno(2.16)$$
Moreover,  appealing to the differential equation (1.11),  we find the equality
$$\nu^2 M_{\nu} (2x) - x^2{d^2\over dx^2}  [M_{\nu} (2x)]  -  x{d\over dx}  [M_{\nu} (2x)] = - 4x^2 [M_{\nu} (2x)] -  \nu e^{-2x}.$$
Hence  the relation (2.16) can be rewritten in the form 
$$ x^{2} {d^2\over dx^2} [ M_{\nu+1} (2x) ]+   \bigg(3  x^2 + 2x(\nu+1) \bigg) {d\over dx} [ M_{\nu+1} (2x) ]+  \bigg(2x^2+ (4\nu+3) x+ \nu(\nu+1) \bigg) M_{\nu+1} (2x) $$

$$ =  - \bigg({x\over 2} - 1\bigg) x^2  {d^2\over dx^2}  [M_\nu(2x) ] - \bigg( x^2+ {5+\nu\over 2} x-1 \bigg) x {d\over dx}  \left[M_\nu(2x) \right] $$

$$-  \bigg(  (\nu+8) x+ 2 \bigg) x M_\nu(2x) -\nu e^{-2x}.\eqno(2.17)$$
On the other hand,  the differential equation (1.11) allows also to deduce the first order differential recurrence relation. Precisely, we have from (2.16), (2.17)

$$ x  {d\over dx} \bigg[  \bigg(3  x + 2\nu+1 \bigg)  M_{\nu+1} (2x) +  \bigg( x^2+ {4+\nu\over 2}  x \bigg)  M_\nu(2x) \bigg]$$

$$ +  \bigg(6x^2+ 4\nu x+ (\nu+1) (2\nu+1)  \bigg) M_{\nu+1} (2x) $$

$$ =   -  \bigg(2x^3+  (\nu+2) x^2 + {\nu (\nu-1)^2\over 2} x- \nu^2\bigg)  M_\nu(2x)-   \left(\nu \left({x\over 2} +1\right) +1\right) e^{-2x}.\eqno(2.18)$$

Further, returning to (2.6) and  invoking (2.11), we write functions $\omega^\alpha_{n,\nu} (2x)$ as follows

$$  \omega^\alpha_{n,\nu} (2x) =  { 2^{\alpha-1} \over \pi} \  \sin(\pi\nu) {d^n\over dx^n}  \bigg[ x^{n+\alpha-1}  e^{-2x} \int_0^\infty  e ^{-t/x} \rho^2_\nu(t)\ t^{-\nu} dt\bigg].\eqno(2.19)$$
Then the differentiation under the integral sign can be performed for $x \ge x_0 >0$ by the absolute and uniform convergence, employing the identity  (see [1], Entry 1.1.3.2)

$${d^n\over dx^n} \left[ x^{\lambda}  e^{-a/x} \right] = (-1)^n n! \  e^{-a/x} x^{\lambda-n} L_n^{-\lambda-1} \left({a\over x}\right).\eqno(2.20)$$
Therefore this gives

$$  \omega^\alpha_{n,\nu} (2x) =  e^{-2x} (2x)^{\alpha+n-1}  n!  \  \sin(\pi\nu)  $$

$$\times   { (-1)^n \over \pi}  \sum_{k=0}^n  {(2x)^{-k}\over (n-k)!}  \int_0^\infty  e ^{-t/x} \rho^2_\nu(t)\  L_k^{-n-\alpha} \left({t\over x}\right) t^{-\nu} dt.\eqno(2.21)$$
For  $\alpha < k-n+1$  the finite sum in (2.21) can be treated via the identity (see [5, Vol. II], Entry 4.4.1.7)

$$  \sum_{k=0}^n  {y^{k} L_k^{-n-\alpha} \left(z\right) \over (n-k)!\ (1-n-\alpha)_k}   = {(y+1)^n\over (1-n-\alpha)_n}\  L_n^{-n-\alpha} \left({yz\over y+1}\right),$$
where $(a)_m$ is the Pochhammer symbol [4]. Precisely, we have
$$\int_0^\infty e^{-2xy} y^{-n-\alpha}  \sum_{k=0}^n  {y^{k} L_k^{-n-\alpha} \left(z\right) \over (n-k)!\ (1-n-\alpha)_k} dy$$

$$= \Gamma (1-n-\alpha) \sum_{k=0}^n  {(2x)^{n+\alpha-k-1}  L_k^{-n-\alpha} \left(z\right) \over (n-k)!} $$

$$=  {\Gamma (1-n-\alpha)\over \Gamma(1-\alpha)} \int_0^\infty e^{-2xy} y^{-n-\alpha} (y+1)^n\  L_n^{-n-\alpha} \left({yz\over y+1}\right) dy.$$
Consequently, comparing with (2.21), we find the following double integral representation

$$  \omega^\alpha_{n,\nu} (2x) =     { (-1)^n  n!  \sin(\pi\nu) \over \pi \Gamma(1-\alpha)}  $$

$$\times  \int_0^\infty  \int_0^\infty  e^{-2x (1+y)- t/x} \rho^2_\nu(t)\ y^{-n-\alpha} (y+1)^n\  L_n^{-n-\alpha} \left({ ty\over x(y+1)}\right)  t^{-\nu} dy dt.$$
On the other hand, recalling the Mellin-Barnes representation (2.9) and differentiating under the integral sign for $x \ge x_0 >0$ due to the same justification, we use the Rodrigues formula (2.4) for the associated Laguerre polynomials and the duplication formula for the gamma function [4] to arrive at the corresponding integral for $ \omega^\alpha_{n,\nu} (2x)$. Indeed, we derive from (2.6), (2.9)

$$ \omega^\alpha_{n,\nu} (2x) =  {e^{-2x}\ 2^\alpha  \sin(\pi\nu)\ n!   \over 4\pi^{2} i} \int_{\gamma-i\infty}^{\gamma+i\infty} \Gamma (s+\nu)\Gamma(s-\nu) \Gamma(1-s)$$

$$\times {\Gamma^2 (s)\over \Gamma(2s)} \    x^{\alpha-s} L_n^{\alpha-s} (2x)  ds.\eqno(2.22)$$
Then making use  the identity (see [5, Vol. II], Entry 4.4.1.1) for associated Laguerre polynomials
$$ L_n^{\beta+1} (x) = \sum_{k=0}^n L_k^\beta (x),$$ 
and replacing $2x$ by $x$,  we easily find from (2.22) the following recurrence relation

$$  \omega^{\alpha+1}_{n,\nu} (x) = x n! \sum_{k=0}^n  {\omega^\alpha_{k,\nu} (x)\over k!}.$$
Furthermore, another identity  for Laguerre polynomials 

$$ L_n^{\beta} (x) =  L_{n}^{\beta+1} (x) -  L_{n-1}^{\beta+1} (x),\eqno(2.23)$$
presumes  the equality

$$  \omega^{\alpha+1}_{n,\nu} (x) =  n  \omega^{\alpha+1}_{n-1,\nu} (x) + x  \omega^{\alpha}_{n,\nu} (x).\eqno(2.24)$$
Recalling (2.21), we employ the closed form for the associated Laguerre polynomials

$$L_n^{\beta} (x)  = \Gamma(1+\beta+n) \sum_{k=0}^n {(-1)^k\over k!} {x^k\over (n-k)! \Gamma(1+\beta+k)},$$
then replace $2x$ by $x$ and use Entries 8.4.23.27, 8.4.46.7 from [5, Vol. III] to represent $ \omega^{\alpha}_{n,\nu} (x) $ in terms of the finite sum, involving integrals with the product of the square of Macdonald's function and Tricomi's function.  Precisely, we obtain

$$\omega^\alpha_{n,\nu} (x) =  {e^{-x}\  x^\alpha \sin(\pi\nu)\ n! \over  \sqrt\pi }  \sum_{k=0}^n {(-1)^k\over k!} {x^k\over (n-k)! } $$

$$\times {1\over 2\pi i} \int_{\gamma-i\infty}^{\gamma+i\infty} {\Gamma (s+\nu)\Gamma(s-\nu)  \Gamma (s) \Gamma(1-s) \Gamma(1+n+\alpha-s)\over \Gamma(s+1/2) \Gamma(1+k+\alpha-s)} \  (2x)^{-s}   ds$$

$$=  {2\over \pi}\   e^{-x}\  x^\alpha \sin(\pi\nu)  \sum_{k=0}^n (-1)^k \binom{n}{k} x^k $$

$$\times \int_0^\infty  K^2_\nu\left(\sqrt{2xt}\right) e^{-t} \Psi\left(k-n,\ 1-\alpha-n; t\right) dt.$$
Meanwhile, the latter Tricomi function can be expressed in terms of the associated Laguerre polynomials via Entry 7.11.4.12 from [5, Vol. III], namely,

$$\Psi\left(k-n,\ 1-\alpha-n; t\right) = (-1)^{n+k} (n-k)! L_{n-k}^{-n-\alpha}(t).$$
Hence we deduce the representation

$$ \omega^\alpha_{n,\nu} (x) =  { 2 (-1)^n\over \pi} e^{-x}\  x^\alpha \sin(\pi\nu)\ n!  \sum_{k=0}^n { x^k \over k!}  \int_0^\infty  K^2_\nu\left(\sqrt{2xt}\right) e^{-t} L_{n-k}^{-n-\alpha}(t) dt,\  |{\rm Re} \nu |  < 1.\eqno(2.25)$$
It is easily seen that this  corresponds to (2.20) via a simple substitution.

On the other hand, employing (1.12), (1.13), (2.4),  we immediately deduce more identities

$$  \omega^{\alpha+1}_{n,\nu+1} (x) -  \omega^{\alpha+1}_{n,\nu-1} (x) =  2\nu  \omega^{\alpha}_{n,\nu} (x) + 2 n! x^\alpha e^{-x} L_n^\alpha (x),$$

$$  \omega^{\alpha+1}_{n,\nu+1} (x) +  \omega^{\alpha+1}_{n,\nu-1} (x) = - 2  \omega^{\alpha}_{n+1,\nu} (x) + 2(n+\alpha +1) \omega^{\alpha}_{n,\nu} (x),$$ 
and, consequently,  

$$  \omega^{\alpha+1}_{n,\nu+1} (x) =  (n+\alpha +\nu+1) \omega^{\alpha}_{n,\nu} (x) -  \omega^{\alpha}_{n+1,\nu} (x) + n! x^\alpha e^{-x} L_n^\alpha (x).$$
Moreover, taking into account (2.24), we have

$$  \omega^{\alpha+1}_{n,\nu+1} (x) =  (\alpha +\nu) \omega^{\alpha}_{n,\nu} (x) -  x\omega^{\alpha-1}_{n+1,\nu} (x) + n! x^\alpha e^{-x} L_n^\alpha (x).$$
Employing the differential equation (1.11), we establish the following recurrence relation

$$  \omega^{\alpha}_{n+2,\nu} (x) -  (2 (n+\alpha)+3)   \omega^{\alpha}_{n+1,\nu} (x)  +   (n+\alpha+1+\nu)   (n+\alpha+1-\nu)  \omega^{\alpha}_{n,\nu} (x) $$

$$ =  \omega^{\alpha+2}_{n,\nu} (x) + \nu n! x^\alpha e^{-x} L_n^\alpha (x).$$
Furthermore, rewriting (2.19) in the form

$$  \omega^\alpha_{n,\nu} (2x) =  { 2^{\alpha-1} \over \pi} \  \sin(\pi\nu) {d^n\over dx^n}  \bigg[ x^{n+\alpha-\nu}  e^{-2x} \int_0^\infty  e ^{-t} \rho^2_\nu(xt)\ t^{-\nu} dt\bigg],\eqno(2.26)$$
we appeal to (2.12), (2.13) and similar arguments as above to deduce from (2.26)

$$ \nu^2 \omega^\alpha_{n,\nu} (2x) =  { 2^{\alpha-1} \over \pi} \  \sin(\pi\nu) {d^n\over dx^n}  \bigg[ x^{n+\alpha-\nu}  e^{-2x} \int_0^\infty  e ^{-t} \bigg[ \rho_{\nu+1}(xt)- xt  \rho_{\nu-1}(xt)\bigg]^2  \ t^{-\nu} dt\bigg]$$

$$=   { 2^{\alpha-1} \over \pi} \  \sin(\pi\nu) {d^n\over dx^n}  \bigg[ x^{n+\alpha-\nu}  e^{-2x} \int_0^\infty  e ^{-t} \rho_{\nu+1}^2 (xt) \ t^{-\nu} dt\bigg]$$

$$-  { 2^{\alpha} \over \pi} \  \sin(\pi\nu) {d^n\over dx^n}  \bigg[ x^{n+\alpha-\nu+1}  e^{-2x} \int_0^\infty  e ^{-t}  \rho_{\nu+1}(xt)\rho_{\nu-1}(xt) \ t^{1-\nu} dt\bigg]$$

$$+  { 2^{\alpha-1} \over \pi} \  \sin(\pi\nu) {d^n\over dx^n}  \bigg[ x^{n+\alpha-\nu+2}  e^{-2x} \int_0^\infty  e ^{-t} \rho_{\nu-1}^2 (xt) \ t^{2-\nu} dt\bigg]$$

$$=   { 2^{\alpha-1} \over \pi} \  \sin(\pi\nu) {d^n\over dx^n}  \bigg[ x^{n+\alpha-\nu+1}  e^{-2x}   {d\over dx} \int_0^\infty  e ^{-t}   \rho^2_{\nu+1}(xt)  \ t^{-\nu} dt\bigg]$$

$$+   2^{\alpha}  {d^n\over dx^n}  \bigg[ x^{n+\alpha+1}  \bigg[ (2x+1) \bigg[ M_{\nu-1}(2x) - M_{\nu+1}(2x) \bigg] \bigg.$$

$$\bigg. + x {d\over dx}  \bigg[M_{\nu-1}(2x) - M_{\nu+1}(2x) \bigg] \bigg]\bigg].$$
Hence (1.12), (2.6)  and again  (2.12) imply

$$ \nu^2 \omega^\alpha_{n,\nu} (2x) = -   2^{\alpha} {d^{n+1}\over dx^{n+1}}  \bigg[  x^{n+\alpha+2}  \bigg[ 2x M_{\nu+1}(2x) +  {d\over dx}  \left[ x M_{\nu+1}(2x) \right]\bigg] \bigg]$$

$$-  2^{\alpha+1}  {d^n\over dx^n}  \bigg[x^{n+\alpha+2}   \bigg[ 2x M_{\nu+1}(2x) +  {d\over dx}  \left[x M_{\nu+1}(2x) \right] \bigg]\bigg] $$

$$+  2^{\alpha} (n+\alpha-\nu+1) {d^n\over dx^n}  \bigg[   x^{n+\alpha+1} \bigg[ 2x M_{\nu+1}(2x) +  {d\over dx}  \left[x M_{\nu+1}(2x) \right]\bigg] $$

$$-   2^{\alpha}  {d^n\over dx^n}  \bigg[ x^{n+\alpha+1}  \bigg[ 2 \left[ e^{-2x} + \nu M_\nu(2x)\right] +  {d\over dx}  \bigg[  e^{-2x} + \nu M_\nu(2x)\bigg] \bigg]\bigg]$$

$$= -  \omega^{\alpha+2}_{n+1,\nu+1} (2x) - {1\over 2} \omega^{\alpha+1}_{n+2,\nu+1} (2x) +  \left(n+\alpha + {3-\nu\over 2}\right)\omega^{\alpha+1}_{n+1,\nu+1} (2x)$$

$$- {1\over 2}  \omega^{\alpha+3}_{n,\nu+1} (2x) + \left(n+\alpha + {3-\nu\over 2}\right)\omega^{\alpha+2}_{n,\nu+1} (2x) $$

$$ -{1\over 2} (n+\alpha-\nu+1) (n+\alpha+1) \omega^{\alpha+1}_{n,\nu+1} (2x)  $$

$$- {\nu\over 2} \omega^{\alpha+1}_{n,\nu} (2x)  -  \nu  \omega^{\alpha}_{n+1,\nu} (2x) +  \nu (n+\alpha+1) \omega^{\alpha}_{n,\nu} (2x).$$
Consequently,  employing  (2.24) and replacing $2x$ by $x$,  we establish the following recurrence relation for functions (2.6)

$$  x \omega^{\alpha}_{n+2,\nu+1} (x) =   \left( 2(\alpha -x) + n+1-\nu\right)\omega^{\alpha+1}_{n+1,\nu+1} (x) -  2\nu x \omega^{\alpha-1}_{n+1,\nu} (x)$$

$$-  \omega^{\alpha+3}_{n,\nu+1} (x) + \left(2\alpha + 1-\nu\right)\omega^{\alpha+2}_{n,\nu+1} (x)  - (n+\alpha-\nu+1) (n+\alpha+1) \omega^{\alpha+1}_{n,\nu+1} (x)  $$

$$- \nu \omega^{\alpha+1}_{n,\nu} (x)   + 2  \nu (\alpha-\nu) \omega^{\alpha}_{n,\nu} (x).$$
Now we recall integral representation (2.11) and appeal to the integral for $\rho_\nu^2$ in terms of $\rho_\nu$ (see [12, formulae (3.11)]. Then, merging  it into the right-hand side of (2.11) and replacing $2x$ by $x$, we  express  $ M_\nu(x)$ as follows

$$  M_\nu(x) =  { e^{-x} \over \sqrt\pi} \  {\sin(\pi\nu)\over \Gamma(\nu+1/2)}  \int_0^\infty  e ^{-t} t^{-\nu} \int_t^\infty \left(y- t\right)^{\nu-1/2}  \rho_\nu\left(2x y \right) {dy dt\over \sqrt y}.\quad - {1\over 2} < \nu < 1.$$
The interchange of the order of integration is guaranteed by Fubini's theorem. The inner integral is calculated in terms of the Kummer function ${}_1F_1$ by virtue of the relations 2.2.2.2 in [2] and 7.11.1.2 in [5, Vol. III]. Thus we derive the representation

$$  M_\nu(x) =  { 2 e^{-x} \over \Gamma(\nu) }  \int_0^\infty   e^{-y}  \rho_\nu (2x y)  {}_1F_1\left({1\over 2} +\nu;\ {3\over 2}; y\right) dy.\eqno(2.27)$$
Moreover, formula (2.20) ((2.25)) can be rewritten accordingly,

$$  \omega^\alpha_{n,\nu} (x) =  2 e^{-x} x^{\alpha}  n!  \  { \sin(\pi\nu) \over \Gamma(\nu+1/2)} $$

$$\times   { (-1)^n \over \sqrt \pi}  \sum_{k=0}^n  {x^{n-k}\over (n-k)!}  \int_0^\infty  e ^{- t} \  L_k^{-n-\alpha} \left( t\right) \int_1^\infty (y^2-1)^{\nu-1/2}  \rho_\nu\left( 2x t y^2\right) dy dt.$$
Taking the integral representation for the scaled Macdonald functions (see [11, formula (2.11)])

$${(-1)^n x^n\over n!} \rho_\nu(x)= \int_0^\infty  u^{\nu+n-1} e^{-u- x/u} L_n^\nu(u) du$$
and the explicit expression for the associated Laguerre polynomials,  we have 

$$  \omega^\alpha_{n,\nu} (x) =  2 e^{-x} x^{\alpha}  n!  \  { \sin(\pi\nu) \over \Gamma(\nu+1/2)} $$

$$\times   { (-1)^n \over \sqrt \pi}  \sum_{k=0}^n \sum_{m=0}^k  {  x^{n-k-m}\  \Gamma( 1+k-n-\alpha) \over 2^m  (n-k)! (k-m)! \Gamma(1+m-n-\alpha)}$$

$$\times \int_0^\infty  e ^{- t}  \int_1^\infty (y^2-1)^{\nu-1/2} y^{-2m} \int_0^\infty  u^{\nu+m-1} e^{-u- (2xty^2)/u} L_m^\nu(u) du dy dt$$

$$=   e^{-x} x^{\alpha}  n!  \  { \sin(\pi\nu) \over \Gamma(\nu+1/2)} $$

$$\times   { (-1)^n \over \sqrt \pi}  \sum_{k=0}^n \sum_{m=0}^k  {  x^{n-k-m}\  \Gamma( 1+k-n-\alpha) \over 2^{m}  (n-k)! (k-m)! \Gamma(1+m-n-\alpha)}$$

$$\times   \int_0^\infty  u^{\nu+m} e^{-u} L_m^\nu(u) \int_1^\infty {(y-1)^{\nu-1/2} y^{-m-1/2} \over u+ 2xy}  dy du,\ {\rm Re}\nu > - {1\over 2},$$
where the interchange of the order of integration is via the dominated convergence theorem. Hence the integral with respect to $y$ is calculated via Entry 2.1.4.7 in [2] in terms of the Gauss hypergeometric function, and we deduce

$$  \omega^\alpha_{n,\nu} (x) =   e^{-x} x^{\alpha}  n!  $$

$$\times   { (-1)^n \over \Gamma(\nu) }  \sum_{k=0}^n \sum_{m=0}^k  {  x^{n-k-m-1}\  \Gamma( 1+k-n-\alpha) (1-\nu)_m \over 2^{m}  (n-k)! (k-m)! (3/2)_m \Gamma(1+m-n-\alpha)}$$

$$\times   \int_0^\infty  u^{\nu+m} e^{-u} L_m^\nu(u)\ {}_2F_1\left(1,\ m+1 -\nu; \  m+ {3\over 2} ;\  - { u\over 2x} \right) du.$$
However, (2.4), integration by parts and the identities    [5, Vol. III]

$${d^n\over dz^n}  \left[ z^n  {}_2F_1\left(1,\ b; c; z \right) \right] =  n!\  {}_2F_1\left(n+1,\ b; c; z \right),$$

$${d^n\over dz^n} \left[ {}_2F_1\left(1,\ b; c; z \right) \right] =  {n! \ (b)_n\over (c)_n} \  {}_2F_1\left(n+1,\ b+n; c+n; z \right)$$
allow to rewrite the latter integral as follows

$$ \int_0^\infty  u^{\nu+m} e^{-u} L_m^\nu(u)\ {}_2F_1\left(1,\ m+1 -\nu; \  m+ {3\over 2} ;\  - { u\over 2x} \right) du$$

$$= {1\over m!}  \int_0^\infty  {d^m\over du^m} \left[ u^{\nu+m} e^{-u} \right] \  u^m {}_2F_1\left(1,\ m+1 -\nu; \  m+ {3\over 2} ;\  - { u\over 2x} \right) du$$

$$= {(-1)^m \over m!}  \int_0^\infty   u^{\nu+m} e^{-u} \   {d^m\over du^m} \left[  u^m {}_2F_1\left(1,\ m+1 -\nu; \  m+ {3\over 2} ;\  - { u\over 2x} \right) \right] du$$

$$= (-1)^m  \int_0^\infty   u^{\nu+m} e^{-u} \   {}_2F_1\left(m+1,\ m+1 -\nu; \  m+ {3\over 2} ;\  - { u\over 2x} \right) du$$

$$= { (2x)^m  (3/2)_m\over m! (1-\nu)_m  }  \int_0^\infty   u^{\nu+m} e^{-u} \  {d^m\over du^m} \left[  {}_2F_1\left(1,\ 1 -\nu; \   {3\over 2} ;\  - { u\over 2x} \right) \right] du$$

$$= { (-1)^m  (2x)^m   (3/2)_m\over (1-\nu)_m }  \int_0^\infty   u^{\nu} e^{-u}  L_m^\nu(u) \  {}_2F_1\left(1,\ 1 -\nu; \   {3\over 2} ;\  - { u\over 2x} \right) du.$$
Consequently, 

$$  \omega^\alpha_{n,\nu} (x) =   e^{-x} x^{\alpha+\nu}  n!  $$

$$\times   { (-1)^n 2^{\nu+1} \over \Gamma(\nu) }  \sum_{k=0}^n \sum_{m=0}^k  { (-1)^k  x^{n-k}\   \Gamma(n-m+\alpha) \over   (n-k)! \ (k-m)!\  \Gamma(n-k+\alpha)}$$

$$\times  \int_0^\infty   u^{\nu} e^{-2xu}  L_m^\nu(2xu) \  {}_2F_1\left(1,\ 1 -\nu; \   {3\over 2} ;\  -  u \right) du.$$
Hence Entries 8.4.49.13, 8.4.33.3, 8.4.23.27 in [5, Vol. III],  the differentiation under the integral and (2.10) sign yield

$$ \int_0^\infty   u^{\nu} e^{-2xu}  L_m^\nu(2xu) \  {}_2F_1\left(1,\ 1 -\nu; \   {3\over 2} ;\  -  u \right) du$$

$$= {(2x)^{-\nu}  \over 4\sqrt \pi i\  m! \Gamma(1-\nu)}  \int_{\gamma-i\infty}^{\gamma+i\infty} {\Gamma (1-s)\Gamma(s) \Gamma(s-\nu)\Gamma(s+\nu)\over \Gamma(s+1/2)} (1-s)_m (2x)^{-s} ds$$

$$= {(2x)^{-\nu}  2^{-m}\over 4\sqrt \pi i\  m! \Gamma(1-\nu)}  {d^m\over dx^m}\int_{\gamma-i\infty}^{\gamma+i\infty} {\Gamma (1-s)\Gamma(s) \Gamma(s-\nu)\Gamma(s+\nu)\over \Gamma(s+1/2)} (2x)^{m-s} ds$$

$$= {(2x)^{-\nu} \over  m! \Gamma(1-\nu)}  {d^m\over dx^m}  \left[ x^m \int_0^\infty e^{-t} K_\nu^2\left(\sqrt{2xt} \right) dt\right]= {(2x)^{-\nu} \Gamma(\nu) \over 2 m! }  {d^m\over dx^m}  \left[ x^m e^x M_\nu(x) \right].$$
Thus we establish the following expansion for the kernels $\omega^\alpha_{n,\nu}$

$$  \omega^\alpha_{n,\nu} (x) =   e^{-x} x^{\alpha} $$

$$\times   \sum_{k=0}^n  (-1)^{n+k}  \binom{n}{k}  {x^{n-k}\over (\alpha)_{n-k}} \sum_{m=0}^k  \binom{k}{m} (\alpha)_{n-m}  {d^m\over dx^m}  \left[ x^m e^x M_\nu(x) \right].$$

$$=   e^{-x} x^{\alpha}    \sum_{m=0}^n   \binom{n}{m} (\alpha)_{n-m}   {d^m\over dx^m}  \left[ x^m e^x M_\nu(x) \right]  \sum_{k=0}^{n-m}  (-1)^{k}  \binom{n-m}{k}  {x^{k}\over (\alpha)_{k}}.$$
But since

$$ (\alpha)_{n-m}  \sum_{k=0}^{n-m}  (-1)^{k}  \binom{n-m}{k}  {x^{k}\over (\alpha)_{k}} = (n-m)!\  L_{n-m}^{\alpha-1} (x),$$
we end up  with the formula

$$  \omega^\alpha_{n,\nu} (x) =  e^{-x} x^{\alpha}  n!   \sum_{m=0}^n {L_{n-m}^{\alpha-1} (x)\over m!}   {d^m\over dx^m}  \left[ x^m e^x M_\nu(x) \right].$$
Moreover, it can be written in terms of the duble sum as follows

$$ \omega^\alpha_{n,\nu} (x) =   x^{\alpha}  n!   \sum_{m=0}^n   \sum_{q=0}^m   {L_{n-m}^{\alpha-1} (x)\over  q! (m-q)!}   {d^q\over dx^q}  \left[ x^m  M_\nu(x) \right]$$

$$ =   x^{\alpha}  n!   \sum_{m=0}^n   \sum_{q=0}^m  L_{n-m}^{\alpha-1} (x)  {\omega^{m-q}_{q,\nu} (x)\over q! (m-q)!}.\eqno(2.28)$$

Further, returning to (1.10) and using simple substitutions, we write the corresponding integral as follows

$$ 2 M_{\nu}(2\sqrt x)= \int_{0}^\infty e^{- \sqrt x ( e^{u}+ e^{-u}) } \left( e^{\nu u} - e^{- \nu u} \right) du $$

$$=   \int_{1}^\infty e^{- \sqrt x ( t  +  t^{-1})} \left( t^{\nu-1} - t^{- \nu -1} \right) dt =   x^{-\nu/2}  \int_{\sqrt x}^\infty e^{- t  - {x\over t} }  t^{\nu-1} dt $$

$$- x^{\nu/2}  \int_{\sqrt x}^\infty e^{- t  - {x\over t} }  t^{-\nu-1} dt.\eqno(2.29)$$
Meanwhile, taking into account the integral representation of the scaled Macdonald functions [11]

$$\rho_\nu(x)=   \int_{0}^\infty e^{- t  - {x\over t} }  t^{\nu-1} dt,\quad x >0\eqno(2.30)$$
and the property $x^\nu \rho_{-\nu}(x) =  \rho_{\nu}(x)$,  we easily find from (2.29)

$$ 2 x^{\nu/2}  M_{\nu}(2\sqrt x) =    \int_{\sqrt x}^\infty e^{- t  - {x\over t} }  t^{\nu-1} dt  - x^{\nu}  \int_{\sqrt x}^\infty e^{- t  - {x\over t} }  t^{-\nu-1} dt$$

$$=  x^{\nu}  \int_0^{\sqrt x} e^{- t  - {x\over t} }  t^{-\nu-1} dt -  \int_0^{\sqrt x} e^{- t  - {x\over t} }  t^{\nu-1} dt$$

$$=  \int_{\sqrt x}^\infty  e^{- t  - {x\over t} }  t^{\nu-1} dt -  \int_0^{\sqrt x} e^{- t  - {x\over t} }  t^{\nu-1} dt = \rho_\nu(x) -   2 \int_0^{\sqrt x} e^{- t  - {x\over t} }  t^{\nu-1} dt.$$
Hence, denoting by $\hat{\rho}_\nu(x)=  2 x^{\nu/2}  M_{\nu}(2\sqrt x) $, we have the following relation

$$\hat{\rho}_\nu(x)= \rho_\nu(x) -   2 x^{\nu/2}  \int_0^{1} e^{-  \sqrt x \left(t  +  {1\over t}\right) }  t^{\nu-1} dt.\eqno(2.31)$$
Invoking (2.13),  we write, using integration by parts

$$\hat{\rho}_{\nu+1}(x)= \rho_{\nu+1}(x) -   2 x^{(\nu+1)/2}  \int_0^{1} e^{-  \sqrt x \left(t  +  {1\over t}\right) }  t^{\nu} dt$$

$$= \nu \rho_{\nu}(x) + x \rho_{\nu-1}(x) +  2 x^{\nu/2} e^{-  2\sqrt x}  -   2 x^{(\nu+1)/2}  \int_0^{1} e^{-  \sqrt x \left(t  +  {1\over t}\right) }  t^{\nu-2} dt$$

$$-  2\nu x^{\nu/2}  \int_0^{1} e^{-  \sqrt x \left(t  +  {1\over t}\right) }  t^{\nu-1} dt =  \nu \hat{\rho}_{\nu}(x) + x \hat{\rho}_{\nu-1}(x) + 2 x^{\nu/2} e^{-  2\sqrt x} ,$$
i.e. we derive an analog of (2.13) for the scaled Macdonald-type functions $\hat{\rho}_{\nu}(x) = 2 x^{\nu/2}  M_{\nu}(2\sqrt x) $

$$ \hat{\rho}_{\nu+1}(x)=  \nu \hat{\rho}_{\nu}(x) + x \hat{\rho}_{\nu-1}(x) + 2 x^{\nu/2} e^{-  2\sqrt x},\ x >0,\ \nu \in \mathbb{C}.\eqno(2.32)$$
{\bf Remark 1}.  Formula (2.31) can be obtained directly from (1.12), replacing $x$ by $2\sqrt x$ and multiplying both sides by $x^{\nu/2}$.

For the derivative we have, accordingly, from (2.13), (2.31), (2.32)

$$x {d\over dx} \left[ \hat{\rho}_\nu(x) \right] = - x\rho_{\nu-1} (x) -   \nu x^{\nu/2}  \int_0^{1} e^{-  \sqrt x \left(t  +  {1\over t}\right) }  t^{\nu-1} dt$$

$$+  x^{(\nu+1)/2}  \int_0^{1} e^{-  \sqrt x \left(t  +  {1\over t}\right) }  t^{\nu} dt +  x^{(\nu+1)/2}  \int_0^{1} e^{-  \sqrt x \left(t  +  {1\over t}\right) }  t^{\nu-2} dt$$

$$= - {x\over 2} \hat{\rho}_{\nu-1} (x)  + {\nu\over 2} \hat{\rho}_{\nu} (x)   -  {1\over 2} \hat{\rho}_{\nu+1} (x) $$

$$=  \nu \hat{\rho}_{\nu} (x)   -   \hat{\rho}_{\nu+1} (x) + x^{\nu/2} e^{-  2\sqrt x} = - x \hat{\rho}_{\nu-1}(x)- x^{\nu/2} e^{-  2\sqrt x} ,$$
i.e. (compare with the formula $\rho^\prime_{\nu} (x) = - \rho_{\nu-1} (x) $) 

$$  {d\over dx} \left[ \hat{\rho}_\nu(x) \right] =   -   \hat{\rho}_{\nu-1} (x) - x^{\nu/2-1} e^{-  2\sqrt x}.\eqno(2.33)$$
Moreover,   identities (2.32), (2.33) immediately imply the relation

$$ \hat{\rho}_{\nu+1}(x) + 2x  {d\over dx} \left[ \hat{\rho}_\nu(x) \right]  =  \nu \hat{\rho}_{\nu}(x) - x \hat{\rho}_{\nu-1}(x),\ x >0,\ \nu \in \mathbb{C},\eqno(2.34)$$
which can be obtained  from (1.13) in a straightforward way, making the corresponding substitutions.    Hence, summing (2.32) and (2.34),  we get, in turn, 

$$ \hat{\rho}_{\nu+1}(x) + x  {d\over dx} \left[ \hat{\rho}_\nu(x) \right]  = \nu \hat{\rho}_{\nu}(x)  + x^{\nu/2} e^{-  2\sqrt x},\ x >0,\ \nu \in \mathbb{C}.\eqno(2.35)$$
Furthermore, we deduce the second order differential recurrence equation for the scaled Macdonald-type functions $\hat{\rho}_{\nu}(x)$. Precisely,  we have from the definition of $ \hat{\rho}_\nu$ and  (2.33)

$$   {d\over dx} \left[ x^{-\nu} \hat{\rho}_\nu(x) \right] = - {d\over dx} \left[  \hat{\rho}_{-\nu}(x) \right] $$

$$=    \hat{\rho}_{-\nu-1} (x) + x^{-\nu/2-1} e^{-  2\sqrt x}=  x^{-\nu-1} \left[ -  \hat{\rho}_{\nu+1} (x) + x^{\nu/2} e^{-  2\sqrt x}\right].\eqno(2.36)$$
Hence

$${d\over dx} \left[ x^{\nu+1}  {d\over dx} \left[ x^{-\nu} \hat{\rho}_\nu(x) \right] \right] = \hat{\rho}_{\nu} (x) + x^{(\nu-1)/2} e^{-  2\sqrt x} $$

$$+ {\nu\over 2}   x^{\nu/2-1} e^{-  2\sqrt x} - x^{(\nu-1)/2} e^{-  2\sqrt x}=  \hat{\rho}_{\nu} (x) +   {\nu\over 2}   x^{\nu/2-1} e^{-  2\sqrt x},\eqno(2.37) $$
and after simplification  the differential recurrence equation takes the final form

$$ x  {d^2\over dx^2} \left[  \hat{\rho}_{\nu}(x) \right] + \left(1-  {\nu\over 2}\right) {d\over dx} \left[  \hat{\rho}_{\nu}(x) \right]-   \hat{\rho}_{\nu}(x) +  {\nu\over 2} \hat{\rho}_{\nu-1}(x) = 0.\eqno(2.38)$$
We note immediate consequences of the equalities (2.33), (2.34), (2.36).  Namely,  we obtain the following differential relations

$$ x^{\nu+1} {d\over dx} \left[ x^{-\nu} \hat{\rho}_\nu(x) \right] + {d\over dx} \left[ \hat{\rho}_{\nu+2}(x) \right] = -2 \hat{\rho}_{\nu+1} (x),\eqno(2.39)$$

$$ x^{\nu+2} {d\over dx} \left[ x^{-\nu-1} \hat{\rho}_{\nu+1}(x) \right] +  x {d\over dx} \left[ \hat{\rho}_{\nu+1}(x) \right] = - \hat{\rho}_{\nu+2} (x) - x  \hat{\rho}_\nu(x),\eqno(2.40)$$

$$ \hat{\rho}_{\nu+2}(x)+ x \hat{\rho}_{\nu}(x) =   (\nu+1) \hat{\rho}_{\nu+1}(x) - 2x  {d\over dx} \left[ \hat{\rho}_{\nu+1}(x) \right].\eqno(2.41)$$

Finally, let us seek a nonnegative twice continuously differentiable function $f_\nu(x),\  \nu \ge 0, x >0$ with the existence of values $f_\nu(0+), f_\nu^\prime(0+)$ for some $\nu >0$  that  satisfies the differential recurrence relations

$$ x^{\nu+1}  { d\over dx} \left[ x^{-\nu} \left[ \hat{\rho}_\nu(x) + f_\nu(x) \right] \right] = - \left[ \hat{\rho}_{\nu+1}(x) + f_{\nu+1}(x) \right],\eqno(2.42)$$

$${d\over dx} \left[ x^{\nu+1}  {d\over dx} \left[ x^{-\nu}\left[  \hat{\rho}_\nu(x) + f_\nu(x) \right] \right] \right] = \hat{\rho}_{\nu} (x) + f_\nu(x).\eqno(2.43) $$
Then  (2.33), (2.36), (2.37), (2.42), (2.43)  imply

$${d\over dx} \left[ x^{\nu+1}  {d\over dx} \left[ x^{-\nu}  f_\nu(x) \right] \right] =  -  {d\over dx} \left[  x^{\nu/2} e^{-  2\sqrt x} + f_{\nu+1}(x) \right],$$

$$ {d\over dx} \left[  f_{\nu+1}(x) \right] =   x^{(\nu-1)/2} e^{-  2\sqrt x} - f_{\nu}(x),$$

$$  x^{\nu+1}  { d\over dx} \left[ x^{-\nu}  f_\nu(x)  \right] = -  f_{\nu+1}(x) - x^{\nu/2} e^{-  2\sqrt x},$$

$${d\over dx} \left[ x^{\nu+1}  {d\over dx} \left[ x^{-\nu}  f_\nu(x) \right]  \right] =  f_\nu(x) -  {\nu\over 2}   x^{\nu/2-1} e^{-  2\sqrt x}.$$
This means that $f_\nu$ obeys the following second order non-homogeneous differential equation

$$x {d^2\over dx^2} \left[ f_\nu(x) \right] + (1-\nu)  {d\over dx} \left[ f_\nu(x) \right] -  f_\nu(x) = -  {\nu\over 2}   x^{\nu/2-1} e^{-  2\sqrt x},\ x >0.\eqno(2.44)$$
When $\nu=0$,  the solution is immediate via  properties of the scaled Macdonald functions [11] and, in particular, the recurrence relation (2.13). So, we find $f_0(x)=\rho_0(x)= 2K_0(2\sqrt x)$.  Generally, it can be solved  by standard methods via Laplace transform, reducing (2.43) to the first order non-homogeneous differential equation. In fact, assuming the existence of the Laplace transform $F_\nu$  of $f_\nu$ 
$$F_\nu(s)= \int_0^\infty e^{-st} f_\nu(t) dt,\quad s >0$$
and its derivatives, we take it on both sides of the equation (2.44) and use integration by parts and differentiation under the integral sign by virtue of the dominated convergence.  As a result we derive

$$ s^2  {d\over ds} \left[  F_\nu(s)  \right] + \left[ (\nu+1) s +1\right]  F_\nu(s)  =  {\nu\over 2} \int_0^\infty e^{-st -  2\sqrt t}\    t^{\nu/2-1}  dt+ (1-\nu)f_\nu(0+).\eqno(2.45)$$
The integral on the right-hand side can be calculated via Entry 2.2.1.6 in [2]  in terms of the parabolic cylinder function $D_\lambda(z)$, and  (2.45)  becomes

 $$   {d\over ds} \left[  F_\nu(s)  \right] + {(\nu+1) s +1\over s^2}\  F_\nu(s)  = {\Gamma(\nu+1)\over 2^{\nu/2} s^{\nu/2+2} } e^{1/ (2s)} D_{-\nu} \left(\sqrt{{2\over s}}\right)+ \left( 1-\nu\right)\ {f_\nu(0+)\over s^2}.$$
It's unique solution is, accordingly,

$$F_\nu(s) = s^{-\nu-1} e^{1/ s} \bigg[  {\Gamma(\nu+1)\over 2^{\nu/2} }  \int_0^s   t^{\nu/2-1} e^{- 1/ (2t)} D_{-\nu} \left(\sqrt{{2\over t}}\right) dt\bigg.$$

$$\bigg.  + \left(1- \nu\right)   f_\nu(0+)\int_0^s t^{\nu-1} e^{-1/t} dt \bigg],\ \nu, s >0.$$
The latter integral is equal to the complementary  incomplete gamma function $\Gamma(-\nu, 1/s)$ (see Entry 2.2.1.11 in [2]). Moreover, in terms of the Whittaker function [5, Vol. III]  the latter equality  reads

$$F_\nu(s) =  {\Gamma(\nu+1)\over 2^{\nu}  s^{\nu+1}}  e^{1/ s} \int_{1/s}^\infty   t^{-\nu/2-5/4} e^{- t/ 2}\  W_{1/4-\nu/2, 1/4} \left( t\right) dt$$

$$+  \left(1- \nu\right)   f_\nu(0+) s^{-\nu-1} e^{1/ s}\   \Gamma\left(-\nu, \ {1\over s}\right).\eqno(2.46)$$
Meanwhile, Entry 2.21.2.2 in [5, Vol. III] suggests the following integral representation of the Whittaker function in terms of the Gauss hypergeometric function
$$W_{1/4-\nu/2, 1/4} \left( t\right)  = {t^{c+1/4-\nu/2}\over \Gamma(c)} e^{-t/2} \int_0^\infty y^{c-1} e^{-ty} {}_2F_1\left({\nu+1\over 2},\  {\nu\over 2};  c;\  -  y \right)dy,$$
where $c >0$ is arbitrary. So, letting $c= \nu+1$, we have the iterated integral from (2.46)

$$F_\nu(s) =     { e^{1/ s} \over s^{\nu+1}}\int_{1/s}^\infty  e^{- t} \bigg[\int_0^\infty \left(y\over 2\right)^{\nu} e^{-ty} {}_2F_1\left({\nu+1\over 2},\  {\nu\over 2}; \  \nu+1;\  -  y \right) dy+  {\left(1- \nu\right)   f_\nu(0+)\over t^{\nu+1}}\bigg] dt.$$
Interchanging the order of integration which is easy to justify by the dominated convergence and calculated the inner integral with respect to $t$, we end up with the equality
$$F_\nu(s) =    \int_0^\infty { e^{-y} \over y s+1} \bigg[ \left(y\over 2\right)^{\nu}   {}_2F_1\left({\nu+1\over 2},\  {\nu\over 2}; \  \nu+1;\  -  y s \right) $$

$$+  {\left(1- \nu\right)   f_\nu(0+)\over (sy +1)^{\nu}}\bigg] dy,\ \nu, \ s >0.\eqno(2.47)$$
To invert the Laplace transform in (2.47), we appeal to Entry 3.35.1.1 in [5, Vol. V] to write the following Laplace transform of its convolution

$$ { 1\over  s+1/y}   {}_2F_1\left({\nu+1\over 2},\  {\nu\over 2}; \  \nu+1;\  -  y s \right) = {  2^{\nu-1}\ \nu \over \sqrt\pi}\    y^{-\nu/2+1/4}$$

$$\times \int_0^\infty e^{- (s+ 1/y) t}   \int_0^t  u^{\nu/2- 5/4}  e^{ u/ (2y) }\  W_{-1/4-\nu/2, 1/4} \left( {u\over y}\right) du.$$
Therefore, substituting the right-hand side of the latter equality in (2.46), interchanging the order of integration and employing the injectivity of the Laplace transform,  we find the desired function $f_\nu(t)$ in the form

$$f_\nu(t) =   {\nu \over 2\sqrt\pi} \int_0^\infty y^{\nu/2- 3/4}  e^{-y- t/y}  \int_0^t  u^{\nu/2- 5/4}  e^{ u/ (2y) }\  W_{-1/4-\nu/2, 1/4} \left( {u\over y}\right) du dy$$

$$+ \left(1- \nu\right)   f_\nu(0+)\ t^\nu  \int_0^\infty  e^{-y-t/y} y^{-\nu-1} dy.$$
But Entry 2.2.2.8 in [2] and (2.30) allow to rewrite the solution above as follows

$$f_\nu(t) =   {\nu \over 2\sqrt\pi \ \Gamma(1+\nu/2)} \int_0^\infty y^{\nu/2- 3/2}  e^{-y- t/y}  \int_0^t  u^{\nu/2- 1/2} \int_0^\infty  { e^{-\xi u/y}\  \xi^{\nu/2} \over (1+\xi)^{(1+\nu)/2} }\  d \xi du dy$$

$$+ \left(1- \nu\right)   f_\nu(0+)\rho_\nu(t),\quad \nu > 0.$$
Interchanging  the order of integration by Fubini's theorem and calculating the integral by $y$,  recalling (2.30), we have after simplification 

$$f_\nu(t) =   {1 \over \sqrt\pi \ \Gamma(\nu/2)} \int_0^\infty    \xi^{(\nu-1)/2} \ \rho_{(\nu-1)/2}\left( t+\xi \right)   \int_0^{t/\xi} {  u^{\nu/2-1} \over (u+ 1)^{(1+\nu)/2} }\  du d \xi$$

$$+ \left(1- \nu\right)   f_\nu(0+)\rho_\nu(t).\eqno(2.48)$$
The integral by $u$ is calculated in [5, Vol. III, Entry 2.1.4.2]. Consequently,  this gives the solution in the form

$$f_\nu(t) =   {t^{\nu/2}  \over \sqrt\pi \ \Gamma(1+\nu/2)} \int_0^\infty   \rho_{(\nu-1)/2}\left( t+\xi \right)  \  {}_2F_1\left({\nu+1\over 2},\  {\nu\over 2}; \  {\nu\over 2} +1;\   -  {t\over \xi}  \right)  {d \xi \over \sqrt \xi}$$

$$+ \left(1- \nu\right)   f_\nu(0+)\rho_\nu(t),\quad \nu,\ t  > 0.$$
Returning to (2.48), we pass through to the limit when $t\to 0+$. Then,  since (see (2.30)) $\rho_\nu(0+) = \Gamma(\nu),\ \nu >0$ and 

$$\int_0^\infty    \xi^{(\nu-1)/2} \ \rho_{(\nu-1)/2}\left( t+\xi \right)   \int_0^{t/\xi} {  u^{\nu/2-1} \over (u+ 1)^{(1+\nu)/2} }\  du d \xi$$

$$\le { 2 t^{\nu/2}\over\nu}   \int_0^\infty  \ \rho_{(\nu-1)/2}\left( \xi \right)   {d\xi\over \sqrt\xi}, $$
where  the latter integral converges for any $\nu > 0$ by virtue of (1.8), (1.9) and the definition of $\rho_\nu(x)= 2 x^{\nu/2} K_\nu(2\sqrt x)$, we arrive at the following algebraic equation

$$  f_\nu(0+) \bigg[1 -(1-\nu)\Gamma(\nu)\bigg] =0.$$
Hence $ f_\nu(0+) =0$ surely for $\nu >1$ and $0 < \nu \le 1, \ \nu\neq \nu_0 \in (0,1)$, where $\nu_0$ is the root of the equation

$$  (1-\nu)\Gamma(\nu) = 1.\eqno(2.49)$$
Consequently, for these values of $\nu$

$$f_\nu(t) =   {t^{\nu/2}  \over \sqrt\pi \ \Gamma(1+\nu/2)} \int_0^\infty   \rho_{(\nu-1)/2}\left( t+\xi \right)  \  {}_2F_1\left({\nu+1\over 2},\  {\nu\over 2}; \  {\nu\over 2} +1;\   -  {t\over \xi}  \right)  {d \xi \over \sqrt \xi}.\eqno(2.50)$$
For instance, let $\nu=1$.  Then (2.48),  Entry 2.16.3.6 in [5, Vol. II] and integration by parts  imply

$$f_1(t) =   {2 \over \pi } \int_0^\infty   \rho_{0}\left( t+\xi \right)  \hbox{arctan}\left(\sqrt{t\over \xi}\right) d \xi$$

$$ =  \rho_1(t) -  {2 \sqrt t \over \pi } \int_t^\infty   {\rho_{1}\left( \xi \right) \over \xi \sqrt{\xi-t}} d \xi =  \rho_1(t) -  {4 \sqrt t \over \pi } \int_{\sqrt t}^\infty   {K_1(2\xi) \over  \sqrt{\xi^2-t}} d \xi =  \rho_1(t) - e^{-2\sqrt t},$$
i.e. $f_1(t) = \rho_1(t) - e^{-2\sqrt t}.$  On the other hand, one can obtain the Mellin-Barnes representation for $f_\nu(t )$ and,  as a consequence, the series expansion via an infinite sum of residues of gamma functions at  simple right-hand poles.  In fact, employing (2.30) and Entry 8.4.49.13 in [5, Vol. III], we substitute it in (2.50),  interchange the order of integration by Fubini's theorem and make simple change of variables. Hence we  derive the equalities

$$f_\nu(t) =   {  t^{\nu} \over \sqrt\pi\  \Gamma(\nu) \ 2\pi i} \int_{\gamma-i\infty}^{\gamma+i\infty}   \rho_{s-\nu}(t)\  \Gamma\left(s+ {1\over 2}\right) \Gamma(\nu- 2s) \   {2^{2s}\over s}  ds$$

$$={  (2t)^{\nu} \over  \Gamma(\nu) \ 4\pi^2 i}\int_0^\infty u^{-\nu-1} e^{-u- t/u}  \int_{\gamma-i\infty}^{\gamma+i\infty}   \Gamma\left(s+ {1\over 2}\right)  \Gamma\left( {\nu+1\over 2}- s\right) \Gamma\left({\nu\over 2} - s\right) \   {u^{s}\over s} ds du$$

$$={  (2t)^{\nu} \over  2 \Gamma(\nu) }\int_0^\infty u^{-\nu/2-1} e^{-u- t/u}  \bigg[ \sum_{n=0}^\infty  { \Gamma\left(n+ (\nu+1)/2\right)  u^{n}  \over  n! (n+\nu/2) \Gamma (1/2+n)} \bigg.$$

$$\bigg. - \sum_{n=0}^\infty  { \Gamma\left(n+ 1+ \nu/2\right)  u^{1/2+n} \over  n! (n+(\nu+1)/2) \Gamma (3/2+n)}  \bigg] du$$

$$=   {(2t)^{\nu} \over  2 \Gamma(\nu) }  \sum_{n=0}^\infty  {1\over n!} \bigg[ { \Gamma\left(n+ (\nu+1)/2\right)  \rho_{n-\nu/2}(t)  \over  (n+\nu/2) \Gamma (1/2+n)}  -  { \Gamma\left(n+ 1+ \nu/2\right)  \rho_{n+ (1-\nu)/2}(t)   \over  (n+(\nu+1)/2) \Gamma (3/2+n)}  \bigg]$$

$$=   {(2t)^{\nu} \over  \Gamma(\nu) \sqrt \pi }  \sum_{n=0}^\infty  {2^{2n-1} \over \Gamma(2n+1) } \bigg[ { \Gamma\left(n+ (\nu+1)/2\right)  \rho_{n-\nu/2}(t)  \over  n+\nu/2 }  -  { \Gamma\left(n+ 1+ \nu/2\right)  \rho_{n+ (1-\nu)/2}(t)   \over  (n+(\nu+1)/2)  (1/2+n)}  \bigg]$$

$$=   {t^{\nu} \over  \Gamma(\nu) \sqrt \pi }  \sum_{n=0}^\infty  2^{n+\nu} (-1)^n  \Gamma\left({n+ \nu+1\over 2}\right)   {\rho_{(n-\nu)/2}(t)  \over  (n+\nu)  n!}, $$
where $\nu, t  > 0, \ 0 < \gamma < \nu/2$ and the latter interchange of the order of integration and summation is due to the estimate for big $n$

$$2^{n}   \Gamma\left({n+ \nu+1\over 2}\right)   {\rho_{(n-\nu)/2}(t)  \over  (n+\nu)\  n!} $$

$$\le   \Gamma\left({n+ \nu+1\over 2}\right)  \Gamma\left({n- \nu\over 2}\right)  {2^n  \over  (n+\nu)\ n!} = O\left( {1\over n^2} \right),\ n \to \infty.$$
Hence, we get the series representation for $f_\nu(t)$

$$f_\nu(t)=  {t^{\nu} \over  \Gamma(\nu) \sqrt \pi }  \sum_{n=0}^\infty  2^{n+\nu} (-1)^n  \Gamma\left({n+ \nu+1\over 2}\right)   {\rho_{(n-\nu)/2}(t)  \over  (n+\nu)  n!},\ \nu, t > 0.\eqno(2.51)$$
Moreover,  recalling (2.30), we establish  the following uniform bounds $(T_0\ge t \ge t_0 >0)$

$$  f_\nu(t) \le  {T_0^{\nu/2} \ 2^\nu \over  \Gamma(\nu) \sqrt \pi } \bigg[ \sum_{n=0}^{[\nu]}  {\left(2\sqrt T_0\right)^{n} \over  (n+\nu)  n!} \  \Gamma\left({n+ \nu+1\over 2}\right)  \Gamma\left({\nu-n\over 2}\right) \bigg.$$   

$$\bigg. +  T_0^{\nu/2}  \sum_{n=[\nu]+1}^{\infty}  {2^{n} \over  (n+\nu)  n!} \  \Gamma\left({n+ \nu+1\over 2}\right)  \Gamma\left({n-\nu\over 2}\right)\bigg],\ \nu \in \mathbb{R}_+\backslash\{\mathbb{N}\},\eqno(2.52)$$

$$  f_\nu(t) \le    {T_0^{\nu/2} \ 2^{2\nu-1} \Gamma(\nu+1/2) \rho_0(t_0) \over  \Gamma^2 (1+\nu) \sqrt \pi } + {T_0^{\nu/2} \ 2^\nu \over  \Gamma(\nu) \sqrt \pi } \bigg[ \sum_{n=0}^{\nu-1}  {\left(2\sqrt T_0\right)^{n} \over  (n+\nu)  n!} \  \Gamma\left({n+ \nu+1\over 2}\right)  \Gamma\left({\nu-n\over 2}\right) \bigg.$$   

$$\bigg. +  T_0^{\nu/2}  \sum_{n=\nu+1}^{\infty}  {2^{n} \over  (n+\nu)  n!} \  \Gamma\left({n+ \nu+1\over 2}\right)  \Gamma\left({n-\nu\over 2}\right)\bigg],\ \nu \in \mathbb{N}.\eqno(2.53)$$

\section{Multiple orthogonal polynomials associated with the Macdonald-type function}

In this section we will investigate  the multiple orthogonal polynomials for the system of weights  $\left( \mu_{\nu,c}(x),  \mu_{\nu+1,c}(x)\right),  \ \mu_{\nu,c}(x)= e^{-cx} \left[ \hat{\rho}_\nu(x)+ f_\nu(x)\right],\  \hat{\rho}_\nu(x) =  2 x^{\nu/2}  M_{\nu}(2\sqrt x),  x, c  >0,  \nu \ge 0,\ \nu \neq \nu_0$ (see (2.49)), where $f_\nu(x)$ is defined by (2.51). In fact, it concerns  the so-called type $I$  polynomials $A_{n,m}^{\nu,c}(x),\ B_{n,m}^{\nu,c}(x),\ n,m \in \mathbb{N}_0$ such that $A_{n,m}^{\nu,c}(x)$ is a polynomial of degree at most $n$,  $B_{n,m}^{\nu,c}(x)$ is a polynomial of degree at most $m$ which satisfy the multiple orthogonality conditions

$$\int_0^\infty \left[ A_{n,m}^{\nu,c}(x)  \mu_{\nu,c}(x) +  B_{n,m}^{\nu,c}(x)  \mu_{\nu+1,c}(x)\right] x^k dx =0,\quad  k= 0,1,2,\dots, n+m\eqno(3.1)$$
and the type $II$ monic polynomials $P_{n,m}^{\nu,c}(x), \ n,m \in \mathbb{N}_0,\ \nu \ge 0, \ c > 0$ such that $P_{n,m}^{\nu,c}$ is a polynomial of degree at most $n+m$ which satisfies the multiple orthogonality conditions

$$\int_0^\infty  P_{n,m}^{\nu,c}(x) \mu_{\nu,c}(x) x^{k} dx=0,\quad k= 0,1,2,\dots,n-1,\eqno(3.2)$$

$$\int_0^\infty  P_{n,m}^{\nu,c}(x) \mu_{\nu+1,c}(x)  x^{k} dx=0,\quad k= 0,1,2,\dots,m-1.\eqno(3.3)$$
Owing to Section 2  each weight function $\mu_{\nu,c}(x)$ is positive and integrable over $\mathbb{R}_+$.  Moreover,  invoking (2.42), (2.43), this gives the following recurrence and differential recurrence relations for the weight functions. Precisely, we have

$$x \mu_{\nu,c}(x) = \mu_{\nu+2,c}(x) - (\nu+1) \mu_{\nu+1,c}(x),\eqno(3.4)$$

$$\mu^\prime_{\nu+1,c}(x) = - c \mu_{\nu+1,c}(x) - \mu_{\nu,c}(x),\eqno(3.5)$$

$$ x \mu^\prime_{\nu,c}(x) = \left(\nu - c x\right)  \mu_{\nu,c}(x)  -  \mu_{\nu+1,c}(x).\eqno(3.6)$$
In the sequel the following notations will be used

$$q_{n,m}^{\nu,c} (x)= A_{n,m}^{\nu,c}(x)  \mu_{\nu,c}(x) +  B_{n,m}^{\nu,c}(x)  \mu_{\nu+1,c}(x),\eqno(3.7)$$

$$P_{2n}^{\nu,c}(x) = P_{n,n}^{\nu,c}(x),\quad\quad  P_{2n+1}^{\nu,c}(x) = P_{n+1,n}^{\nu,c}(x).\eqno(3.8)$$

Furthermore,  for the weight moments 

$$d^\nu_{n,c}= \int_0^\infty x^n \mu_{\nu,c}(x) dx,\quad n \in \mathbb{N}_0.\eqno(3.9)$$
using (3.5), (3.6) and making integration by parts,  we derive the recurrence relations

$$d^\nu_{n,c} = n d^{\nu+1}_{n-1,c} - c d^{\nu+1}_{n,c},\eqno(3.10)$$

  $$ (n+1+\nu) d^\nu_{n,c}= d^{\nu+1}_{n,c} + c d^\nu_{n+1,c}.\eqno(3.11)$$
Hence, writing (3.10) for $n-1$, multiplying both sides by $n$  and subtracting (3.9), we get

 $$ n (n+\nu) d^\nu_{n-1,c} =  c d^{\nu+1}_{n,c}  + (c n+1) d^\nu_{n,c}.$$
Appealing again to (3.10), we write the final three term recurrence relation 

$$ c^2 d^\nu_{n+1,c}  =  \left( c (2n+1+\nu) +1\right) d^\nu_{n,c} - n (n+\nu) d^\nu_{n-1,c},\ n \in \mathbb{N}_0,\eqno(3.12)$$
where we put $d^\nu_{-1,c} =0$. We observe from (3.8) that  $d\ [d^\nu_{n,c}] / dc= - d^\nu_{n+1,c}.$  Then, differentiating through  (3.12) with respect to $c$, we find

$$   c^2 d^\nu_{n+2,c}  =  \left( c (2n+3+\nu) +1\right) d^\nu_{n+1,c} - (n+1)(n+1 + \nu )  d^\nu_{n,c}$$
which the same equality when we replace $n$ by $n+1$ in (3.12).  Moreover, denoting by $y(c)= d^\nu_{n,c}$, it is straightforward to verify from the previous equality that $y(x)$ satisfies the following second order differential equation

$$  x^2 y^{\prime\prime} +  \left( x (2n+3+\nu) +1\right) y^{\prime} +  (n+1)(n+1 + \nu )  y = 0.\eqno(3.13)$$
Meanwhile, the moments can be calculated in a straightforward way.  Indeed, recalling (1.10), (2.51) and the definition of the scaled Macdonald functions, we deduce 

$$d_{n,c}= \int_0^\infty x^n \mu_{\nu,c}(x) dx =  4  \int_0^\infty \sinh(\nu u) \int_0^\infty x^{2n+\nu+1} e^{-cx^2 - 2x \cosh u} dx du  $$

$$+   {4\over  \Gamma(\nu) \sqrt \pi }  \sum_{k=0}^\infty  {2^{k+\nu} (-1)^k \over  (k+\nu)  k!}  \Gamma\left({k+ \nu+1\over 2}\right)    \int_0^\infty x^{2n+1+ (k+ 3\nu)/2} e^{-cx^2} K_{(k-\nu)/2}\left(2 x\right) dx,$$
where the interchange  the order of integration and integration and summation  are guaranteed by  the dominated convergence (see (1.1), (1.7), (1.8), (1.9),\  $\nu \ge 0,\  n \in \mathbb{N}_0$, (2.52), (2.53)) as, for instance, 

$$\int_0^\infty \sinh(\nu u) \int_0^\infty x^{n+\nu/2} e^{-cx - 2\sqrt x \cosh u} dx du \le \int_0^\infty x^{n+\nu/2} e^{-cx} K_\nu(  2\sqrt x)  dx < \infty.$$
Appealing  to (2.45), (2.46) an Entry 3.14.3.10 in [2], we get the formulae 

$$ \int_0^\infty x^{2n+\nu+1} e^{-cx^2 - 2 x \cosh u} dx  = {\Gamma(2(n+1)+\nu)\over (2c)^{n+1+\nu/2}}\  e^{\cosh^2 (u)/ 2c} D_{-\nu- 2(n+1)} \left({\sqrt 2 \cosh u\over \sqrt c}\right) $$

$$= {\Gamma(2(n+1)+\nu)\over  2^{2(n+1)+\nu}\  c^{n+\nu/2+3/4}} {e^{\cosh^2 u / 2c}\over \sqrt{\cosh u}} W_{-\nu/2- n -3/4,\ 1/4 } \left({ \cosh^2 u\over  c}\right), $$

$$ \int_0^\infty x^{2n+1+ (k+ 3\nu)/2} e^{-cx^2} K_{(k-\nu)/2}\left(2 x\right) dx =   c^{- (4n+k+ 3\nu+2)/4} \  { e^{1/c}\over 4}\  \Gamma\left(n+1 +\nu\right) $$

$$\times \Gamma\left(n+1+ {k+\nu\over 2}\right) W_{- (4n+k+ 3\nu+2)/4,\  (k-\nu)/4} \left({1\over c}\right). $$
Thus the moments $d_{n,c}$ can be written in the form

$$d_{n,c}=    {\Gamma(2(n+1)+\nu)\over  2^{2n+\nu}\  c^{n+\nu/2+3/4}} \int_0^\infty {\sinh(\nu u)\over \sqrt{\cosh u}} \  e^{\cosh^2 u/ 2c}\  W_{-\nu/2- n -3/4,\ 1/4 } \left({ \cosh^2 u\over  c}\right) du$$

$$+ { e^{1/c}\  \Gamma\left(n+1 +\nu\right) \over  \Gamma(\nu) c^{ (4n+ 3\nu+2)/4} \ \sqrt \pi}  \sum_{k=0}^\infty  {2^{k+\nu} (-1)^k \over c^{k/4}  (k+\nu)  k!}\   \Gamma\left({k+ \nu+1\over 2}\right)  $$

$$\times \Gamma\left(n+1+ {k+\nu\over 2}\right) W_{- (4n+k+ 3\nu+2)/4,\  (k-\nu)/4} \left({1\over c}\right).\eqno(3.14)$$
Meanwhile, in the case $\nu > 0, \nu \notin \mathbb{N}$ via (2.30), (2.31), Entries  2.1.2.3, 2.2.2.8, 3.8.2.7  in [2] and the Mellin-Parseval equality [7] we derive

$$\int_0^\infty x^n e^{-cx} \hat{\rho}_\nu(x) dx = \int_0^\infty x^n e^{-cx}  \rho_\nu(x) dx - 2 \int_0^\infty x^{n+\nu/2}  e^{-cx} \int_0^{1} e^{-  \sqrt x \left(t  +  {1\over t}\right) }  t^{\nu-1} dtx$$

$$ = \int_0^\infty x^n e^{-cx}  {1\over 2\pi i} \int_{\gamma-i\infty}^{\gamma+i\infty} \Gamma(s)\Gamma(s+\nu) x^{-s} ds dx - 2 \sum_{k=0}^\infty {(-1)^k\over k!} \int_0^\infty x^{n+(k+\nu)/2}  e^{-cx-\sqrt x} $$

$$\times \int_0^{1} e^{-  \sqrt x {1-t\over t} }  t^{k+\nu-1} dtdx =  {c^{-n-1} \over 2\pi i} \int_{\gamma-i\infty}^{\gamma+i\infty} \Gamma(s)\Gamma(s+\nu) \Gamma(n+1-s)  c^{s} ds $$

$$- 2 \sum_{k=0}^\infty {(-1)^k\over k!} \int_0^\infty x^{n+(k+\nu)/2}  e^{-cx}  \int_0^{\infty}  {  e^{-  (t+1) \sqrt x} \over (1+t)^{k+\nu+1}}\  dt dx$$

$$=   n! \int_0^\infty {e^{-t}\  t^{n+\nu}   dt\over (ct + 1)^{n+1} } - 4 \sum_{k=0}^\infty {(-1)^k\over k!} \int_0^\infty x^{2(n+k+\nu)+1}  e^{-cx^2}  \int_{ x}^{\infty}  e^{-  t}  t^{-k-\nu-1}\  dt dx$$

$$=   {n! \ \Gamma(n+\nu+1) \over c^{n+(\nu+1)/2} }\ e^{1/(2c)}\  W_{-n- (\nu+1)/2,\  \nu/2} \left({1\over c}\right) $$

$$- {2\over c^{n+1+\nu/2}} \sum_{k=0}^\infty {(-1)^k\over c^{k/2} k!} \bigg[ \pi (-1)^{k+1}   {(k+\nu+1)_n\over  c^{(k+\nu)/2}\  \sin(\pi\nu)}\bigg.$$

$$\bigg.+ {1\over 2} \ (1+(k+\nu)/2)_{n}\  \Gamma\left( {k+\nu\over 2}\right) {}_2F_2\left( - {k+\nu\over 2},\ n+1+{k+\nu\over 2} ; {1\over 2},\  {2-k-\nu\over 2};\ {1\over 4c} \right)$$

$$+ {((k+\nu+1)/2)_{n+1}  \over (1-k-\nu)  \sqrt c} \Gamma\left(  {k+\nu+1\over 2}\right) {}_2F_2\left(  {1-k-\nu\over 2},\ n+{k+\nu+3\over 2} ; {3\over 2},\  {3-k-\nu\over 2};\ {1\over 4c} \right)\bigg],$$
where $\gamma >0$ and interchanges of the order of integration and summation are guaranteed by the dominated convergence theorem.  Hence, combining with (3.9), we get the moments in the form

$$d_{n,c} =   {n! \ \Gamma(n+\nu+1) \over c^{n+(\nu+1)/2} }\ e^{1/(2c)}\  W_{-n- (\nu+1)/2,\  \nu/2} \left({1\over c}\right) $$

$$+ { e^{1/c}\  \Gamma\left(n+1 +\nu\right) \over  c^{ (4n+ 3\nu+2)/4} }  \sum_{k=0}^\infty  { (-1)^k  (\nu)_k \over c^{k/4}\    k!}  (1+(k+\nu)/2)_{n} W_{- (4n+k+ 3\nu+2)/4,\  (k-\nu)/4} \left({1\over c}\right)$$

$$- {2\over c^{n+1+\nu/2}}  \sum_{k=0}^\infty {(-1)^k\over c^{k/2} k!} \bigg[ \Gamma(1-\nu)\  { (-1)^{k+1}\ (\nu)_{n+k+1}  \over  c^{(k+\nu)/2} \ (\nu)_{k+1}   }+ {1\over 2} \ (1+(k+\nu)/2)_{n}\  \Gamma\left( {k+\nu\over 2}\right)$$

$$\times  {}_2F_2\left( - {k+\nu\over 2},\ n+1+{k+\nu\over 2} ; {1\over 2},\  {2-k-\nu\over 2};\ {1\over 4c} \right)+ {((k+\nu+1)/2)_{n+1}  \over (1-k-\nu)  \sqrt c} \Gamma\left(  {k+\nu+1\over 2}\right)$$

$$\times  {}_2F_2\left(  {1-k-\nu\over 2},\ n+{k+\nu+3\over 2} ; {3\over 2},\  {3-k-\nu\over 2};\ {1\over 4c} \right)\bigg].\eqno(3.15)$$

\section{Differential properties}

In this section we establish differential properties for polynomials of the type $I$ and $II$ $q_{n,m}^{\nu,c},\   P_{n,m}^{\nu,c}$, respectively,  when $n=m$ or $n=m+1$.  We begin with

{\bf Theorem 1}. {\it Let $c >0,\ \nu \ge 0,  n \in \mathbb{N}$. For the type $I$ multiple orthogonal polynomials it has the following differential recurrence relations}

$${d\over dx} \left[  q_{n,n}^{\nu+1,c} (x)\right] = q_{n+1,n}^{\nu,c} (x),\quad\quad\quad\quad  {d\over dx}  \left[ q_{n+1,n}^{\nu+1,c} (x)\right]= q_{n+1,n+1}^{\nu,c} (x). \eqno(4.1)$$

\begin {proof} In fact, using (3.4), (3.5), (3.6),  (3.7), we have

$${d\over dx}  \left[q_{n,n}^{\nu+1,c} (x)\right]=  R_{n,n}^{\nu,c}(x)  \mu_{\nu,c}(x) +  S_{n,n}^{\nu,c}(x)  \mu_{\nu+1,c}(x),$$

$${d\over dx}  \left[ q_{n+1,n}^{\nu+1,c} (x)\right] =  R_{n+1,n}^{\nu,c}(x)  \mu_{\nu,c}(x) +  S_{n+1,n}^{\nu,c}(x)  \mu_{\nu+1,c}(x),$$
where

$$R_{n,m}^{\nu,c}(x) = x {d\over dx} B_{n,m}^{\nu+1,c}(x) -  A_{n,m}^{\nu+1,c}(x) - cx B_{n,m}^{\nu+1,c}(x),$$

$$S_{n,m}^{\nu,c}(x)  = {d\over dx} A_{n,m}^{\nu+1,c}(x) + (\nu+1)  {d\over dx} B_{n,m}^{\nu+1,c}(x)  - (c(\nu+1)+1) B_{n,m}^{\nu+1,c}(x) $$

$$-   c A_{n,m}^{\nu+1,c}(x),$$ 
and  $R_{n,n}^{\nu,c}(x), R_{n+1,n}^{\nu,c}(x),  S_{n+1,n}^{\nu,c}(x)$ are polynomials of degree at most $n+1$ and $S_{n,n}^{\nu,c}(x)$ is a polynomial of degree at most $n$. As we see from (3.7),  the expression $x^k q_{n,m}^{\nu+1,c}(x),\ k \in \mathbb{N}_0, \ \nu \ge 0$ vanishes at infinity and near the origin.  Moreover, the integration by parts in (3.1), where we replace $\nu$ by $\nu+1$ and a pair $(n,m)$ by $(n,n), (n+1,n)$, accordingly, shows that for positive integers $n$ the following orthogonality properties take place

$$\int_0^\infty {d\over dx}  \left[q_{n,n}^{\nu+1,c} (x) \right]  x^k dx = 0,\quad k=0,1,\dots, 2n+1,$$

$$\int_0^\infty {d\over dx}  \left[ q_{n+1,n}^{\nu+1,c} (x) \right]  x^k dx = 0,\quad k=0,1,\dots, 2n+2,$$
which guarantee (4.1) and complete the proof.

\end {proof} 

For the type $II$ polynomials we have the following result.

{\bf Theorem 2}. {\it Let $c >0,\ \nu \ge 0,  n \in \mathbb{N}$. For the type $II$ multiple orthogonal polynomials the differential recurrence relation has  the form}

$${d\over dx} \left[ P_{n}^{\nu,c}(x)\right] = n P_{n-1}^{\nu+1,c}(x). \eqno(4.2)$$

\begin {proof} In view of (3.2), (3.3), (3.8), the integration by parts yields

$$\int_0^\infty {d\over dx} \left[ P_{n}^{\nu,c}(x)\right] x^k \mu_{\nu+1,c}(x) dx =  - k  \int_0^\infty P_{n}^{\nu,c}(x) x^{k-1} \mu_{\nu+1,c}(x) dx$$

$$-  \int_0^\infty P_{n}^{\nu,c}(x) x^{k} \mu^\prime_{\nu+1,c}(x) dx,$$
where the integrated terms vanish owing to the behavior of $\mu_{\nu+1,c}$ near $0$ and $\infty$.  Consequently, this gives

$$\int_0^\infty {d\over dx} \left[ P_{n,n}^{\nu,c}(x)\right] x^k \mu_{\nu+1,c}(x) dx = 0,\quad k=0,1,\dots, n-1,\eqno(4.3)$$

$$\int_0^\infty {d\over dx} \left[ P_{n+1,n}^{\nu,c}(x)\right] x^k \mu_{\nu+1,c}(x) dx = 0,\quad k=0,1,\dots, n-1.\eqno(4.4)$$
Analogously, with the use of (3.4), (3.5), we derive

$$\int_0^\infty {d\over dx} \left[ P_{n}^{\nu,c}(x)\right] x^k \mu_{\nu+2,c}(x) dx =  - k  \int_0^\infty P_{n}^{\nu,c}(x) x^{k-1} \mu_{\nu+2,c}(x) dx$$

$$-  \int_0^\infty P_{n}^{\nu,c}(x) x^{k} \mu^\prime_{\nu+2,c}(x) dx$$

$$=   - k  \int_0^\infty P_{n}^{\nu,c}(x) x^{k-1} \left[ x \mu_{\nu,c}(x)  + (\nu+1) \mu_{\nu+1,c}(x) \right] dx$$

$$+ \int_0^\infty P_{n}^{\nu,c}(x) x^{k}\left[ c \mu_{\nu+1,c}(x) + \mu_{\nu,c}(x)\right]  dx.$$
Hence, this implies the orthogonality conditions in the form

$$\int_0^\infty {d\over dx} \left[ P_{n,n}^{\nu,c}(x)\right] x^k \mu_{\nu+2,c}(x) dx = 0,\quad k=0,1,\dots, n-2,$$

$$\int_0^\infty {d\over dx} \left[ P_{n+1,n}^{\nu,c}(x)\right] x^k \mu_{\nu+2,c}(x) dx = 0,\quad k=0,1,\dots, n-1.$$
Meanwhile, we are  considering the monic polynomials.  Thus we establish (4.2) and complete the proof of Theorem 2.

\end{proof}

Now, making the following notations

$$Q_{2n}^{\nu,c}(x) = q_{n,n}^{\nu,c}(x),\quad\quad\quad  Q_{2n+1}^{\nu,c}(x) = q_{n+1,n}^{\nu,c}(x),$$
we establish the Rodrigues formula for the type $I$ polynomials.  Precisely, we have 

{\bf Theorem 3}. {\it For the type $I$ multiple orthogonal polynomials  the following Rodrigues formula takes place

$$Q_{n}^{\nu,c}(x) = {d^{n+1}\over dx^{n+1}}\   \mu_{\nu+n+1,c}(x),\quad\quad \nu \ge 0,\quad c >0,\quad n\in \mathbb{N}_0\eqno(4.5)$$
up to a multiplicative factor which is chosen    to be one. }

\begin{proof} Let $n=0$. Then (3.5) yields $ A_{0,0}^{\nu,c}(x) = -1,\  B_{0,0}^{\nu,c}(x) = -c$ and (4.5) holds.  Assuming (4.5) for $n=m$ and all $\nu \ge 0,\  c >0$, we take $n=m+1$ and appeal to (4.1) to obtain

$$ Q_{m+1}^{\nu,c}(x) = {d\over dx}\  Q_{m}^{\nu+1,c}(x) =  {d^{m+2}\over dx^{m+2}}\   \mu_{\nu+m+2,c}(x) $$
which proves the Rodrigues formula by induction.

\end{proof}

{\bf Corollary 1}. {\it The sequence $\{Q_{n}^{\nu,c}(x)\}_{n\ge 0}$ obeys the following recurrence relation}

$$ Q_{n+1}^{\nu,c}(x) = - c Q_{n}^{\nu+1,c}(x) - Q_{n}^{\nu,c}(x),\quad\quad \nu \ge 0,\quad c >0,\quad n\in \mathbb{N}_0.\eqno(4.6)$$

\begin{proof} In fact, appealing to (4.5) and (3.5) we write in a straightforward way

$$ Q_{n+1}^{\nu,c}(x) =  -  {d^{n+1}\over dx^{n+1}}\  \left( c \mu_{\nu+n+2,c}(x) +  \mu_{\nu+n+1,c}(x) \right) = - c Q_{n}^{\nu+1,c}(x) - Q_{n}^{\nu,c}(x). $$
 
\end{proof}

{\bf Theorem 4}. {\it Let $\nu \ge 0,\  c >0,\quad n\in \mathbb{N}_0$. The sequence $\{Q_{n}^{\nu,c}(x)\}_{n\ge 0}$ has the representation }

$$ Q_{n}^{\nu,c}(x) = (-1)^{n+1} \sum_{k=0}^{n+1} \binom{n+1}{k} c^k \mu_{\nu+k,c}(x).\eqno(4.7)$$

\begin{proof}  The proof follows immediately, employing the fundamental combinatorial identity

$$\binom{n+2}{k} = \binom{n+1}{k} + \binom{n+1}{k-1}$$  
and checking (4.6) for the expression (4.7), where we replace $n$ by $n+1$. 

\end{proof} 

{\bf Corollary 2}. {\it The generating function for polynomials $Q_{n}^{\nu,c}$ has the series representation

$$ H(x,z)= e^{-z} \sum_{k=0}^\infty {(-cz)^k\over k!}\mu_{\nu+k,c}(x),\quad \nu > 0,\  x, c >0,\eqno(4.8)$$
where the series converges in the domain  $|z| < 1/c$. }

\begin{proof} In fact, using (4.7), where $Q_{-1}^{\nu,c}(x) = \mu_{\nu,c}(x)$,  we derive

$$H(x,z)= \sum_{n=0}^\infty   Q_{n-1}^{\nu,c}(x) {z^n\over n!} = \sum_{n=0}^\infty   {(-z)^n\over n!} \sum_{k=0}^{n} \binom{n}{k} c^k \mu_{\nu+k,c}(x) $$

$$=  \sum_{k=0}^\infty  {c^k\over k!}\   \mu_{\nu+k,c}(x) \sum_{n=k}^{\infty}  {(-z)^n\over (n-k)!} =  e^{-z} \sum_{k=0}^\infty {(-cz)^k\over k!}\mu_{\nu+k,c}(x).$$
Meanwhile, from (1.10), (2.31) and (2.52), (2.53)  we find that

$$ \mu_{\nu+k,c}(x) = e^{-cx}\left[ O\left(\Gamma(\nu+k)\right) + O\left( x^{\nu/2} \rho_0(x)\right)\right.$$

$$\left. + O\left( (4x)^{(\nu+k)/2} e^{2\sqrt x}\right) + O\left( (2x)^{\nu+k} \right)\right],\quad k \to \infty.$$
Therefore, the convergence of the latter series is plainly guaranteed in the domain $|z| < 1/c$.

\end{proof} 

 {\bf Remark 2}.    Theorem 4 and  recurrence relation (3.4) give a starting point to establish  explicit formulas for the type I polynomials which we leave for an interested reader.  
 
 {\bf Remark 3}. It is known that the type II multiple orthogonal polynomials satisfies a third order recurrence relation of the form
 
 $$x P_{n}^{\nu,c}(x) = P_{n+1}^{\nu,c}(x) + b_n P_{n}^{\nu,c}(x) + c_n P_{n-1}^{\nu,c}(x) + d_n P_{n-2}^{\nu,c}(x).$$
 Using the formula (3.15) for the moments, it is possible to find explicit formulas for the coefficients  $b_n, c_n, d_n$ as well as an explicit formula for polynomials $P_{n}^{\nu,c}$ and the generating function.  These problems we remain for the engaged reader.

\bigskip
\centerline{{\bf Acknowledgments}}
\bigskip

The work was partially supported by CMUP, which is financed by national funds through FCT (Portugal) under the project with reference UID/00144/2025.

\bigskip
\centerline{{\bf References}}
\bigskip
\baselineskip=12pt
\medskip
\begin{enumerate}

 \item[{\bf 1.}\ ] Yu. A. Brychkov, {\it Handbook of Special Functions: Derivatives, Integrals, Series and Other Formulas}, CRC Press, Chapman and Hall, 2008.

 \item[{\bf 2.}\ ]
Yu. A. Brychkov,  O.I. Marichev,  N.V.  Savischenko, {\it Handbook of Mellin transforms. Advances in Applied Mathematics},  CRC Press, Boca Raton, FL, 2019. 

\item[{\bf 3.}\ ]  M. Krakowski, On certain functions connected with the Bessel functions,  {\it Zastos. Mat.} {\bf 4} (1958), 130-141 (in Polish).

\item[{\bf 4.}\ ]  N.N. Lebedev, {\it Special Functions and Their Applications}, Dover,  New York, 1972.

\item[{\bf 5.}\ ] A.P. Prudnikov, Yu.A. Brychkov and O.I. Marichev, {\it Integrals and Series}. Vol. I: {\it Elementary Functions}, Vol. II: {\it Special Functions}, Gordon and Breach, New York and London, 1986,  Vol. III: {\it More Special Functions}, Gordon and Breach, New York and London, 1990. Vol. V: {\it Inverse Laplace Transforms}, Gordon and Breach, New York and London, 1992.

\item[{\bf 6.}\ ]  S. Yakubovich and L. Gusarevich, On the non-convolution transformation with Macdonald type kernel function,
 {\it Fract. Calculus and Appl. Anal.} {\bf 1} (1998), N 3,  297- 309.

\item[{\bf 7.}\ ]S. Yakubovich and Yu. Luchko, {\it The Hypergeometric Approach to Integral Transforms and Convolutions},
(Kluwers Ser. Math. and Appl.: Vol. 287), Dordrecht, Boston, London, 1994.

\item[{\bf 8.}\ ] S. Yakubovich, {\it Index Transforms}, World Scientific Publishing Company, Singapore, New Jersey, London and Hong Kong, 1996.

\item[{\bf 9.}\ ] S. Yakubovich, A class of integral equations and index transformations related to the modified and incomplete Bessel functions. J. Integral Equations Appl. 22 (2010), no. 1, 141-164.

\item[{\bf 10.}\ ] S. Yakubovich, Certain identities, connection and explicit formulas for the Bernoulli and Euler numbers and the Riemann zeta-values, {\it Analysis} {\bf 35} (2015), 1,  59- 71.

\item[{\bf 11.}\ ] S. Yakubovich, Orthogonal polynomials with ultra-exponential weight functions: an explicit solution to the Ditkin-Prudnikov problem. {\it Constr. Approx. } {\bf 53} (2021), no. 1, 1-38.

\item[{\bf 12.}\ ] S. Yakubovich, A method of composition orthogonality and new sequences of orthogonal polynomials and functions for non-classical weights. {\it J. Math. Anal. Appl. } {\bf 499} (2021), 125032.

\end{enumerate}

\vspace{5mm}

\noindent S.B.Yakubovich\\
Department of  Mathematics,\\
Faculty of Sciences,\\
University of Porto,\\
Campo Alegre st., 687\\
4169-007 Porto\\
Portugal\\
E-Mail: syakubov@fc.up.pt\\

\end{document}